\documentclass{amsart}[12pt]
\usepackage{amsmath, amstext, amssymb, amsopn, amsthm, amscd, amsfonts, times, mathptmx, epsfig}

\newcommand{\D}{\mathbb D}
\newcommand{\R}{\mathbb R}
\newcommand{\C}{\mathbb C}
\newcommand{\HP}{\mathbb H}
\newcommand{\OI}{\mathbb O}
\newcommand{\N}{\mathbb N}
\newcommand{\Q}{\mathbb Q}
\newcommand{\Z}{\mathbb Z}
\newcommand{\SC}{\mathbb S}
\newcommand{\SO}{\widehat{\mathbb S}}
\newcommand{\T}{\mathbb T}

\newtheorem{theo}{Theorem}
\newtheorem{lemm}{Lemma}
\newtheorem{prop}{Proposition}
\newtheorem{coro}{Corollary}
\newtheorem*{assu}{Assumption}
\newtheorem{conj}{Conjecture}

\theoremstyle{definition}
\newtheorem{defi}{Definition}

\newtheorem{note}{Note}

\newtheorem{exam}{Example}

\newcommand{\gzg}[1]{[\negmedspace [ #1 ]\negmedspace ]}

\newcommand{\cal}{\mathcal}
\newcommand{\bast}{{}^{\ast}}
\newcommand{\bstar}{{}^{\text{\tiny $\bullet$}}}

\newcommand{\gcom}[1]{\gzg{\,\widetilde{ #1 }\,}}

\newcommand{\fgrm}{\gzg{\pi}_{1}}

\title{The Algebraic Theory of the Fundamental Germ}
\author{T.M. Gendron}
\address{Instituto de Matem\'{a}ticas, Universidad
Nacional Aut\'{o}noma de M\'{e}xico, Unidad Cuernavaca, Av.
Universidad S/N, C.P. 62210 Cuernavaca, Morelos, M\'{E}XICO}
\email{tim@matcuer.unam.mx} 
\date{25 April 2005}
\begin{document}
\vspace{2cm}
\begin{abstract} This paper introduces a notion of fundamental
group appropriate for laminations.
\end{abstract}
\maketitle

\maketitle
\section*{Introduction}  Let $\cal L$ be a lamination: a
space modeled on a ``deck of cards'' ${\R}^{n}\times {\sf T}$,
where ${\sf T}$ is a topological space and overlap homeomorphisms take cards
to cards continuously in the deck direction $\sf T$. One thinks of
$\cal L$ as a family of manifolds, the leaves, bound by a transversal
topology prescribed locally by ${\sf T}$. Using this picture, many
constructions familiar to the theory of manifolds
can be extended to laminations via the ansatz:

\begin{quote}
{\small Replace manifold object $A$ by a family of manifold
objects $\{ A_{L}\}$ existing on the leaves of $\cal L$ and
respecting the transverse topology.}
\end{quote}

For example, one defines a smooth structure to be a family of
smooth structures on the leaves in which the card gluing
homeomorphisms occurring in a pair of overlapping decks vary
transversally in the smooth topology. Continuing in this way, 
constructions over $\R$, such as tensors, de Rham
cohomology groups, {\em etc.} may be defined.

Identifying those constructions classically defined over $\Z$ is
not as straightforward, especially if one wishes to follow tradition
and define them geometrically.  To see why this is true, 
consider the case of an exceptionally well-behaved lamination:
an inverse limit 
$\widehat{M}=\lim_{\longleftarrow}M_{\alpha}$ of manifolds by covering maps.
Such a system induces a direct
limit of de Rham cohomology groups, and there is a canonical map
from this limit into the tangential cohomology groups
$H^{\ast}(\widehat{M};\; {\R})$ with dense image. In fact,
here one may use the system to define -- by completion of limits --
tangential homology groups $H_{\ast}(\widehat{M};\;\R )$ as well. If one
endeavors to use this point of view to define the groups $\pi_{1}$,
$H_{\ast}(\cdot \;\;\Z )$, $H^{\ast}(\cdot \;\;\Z )$, the result is
failure since the systems they induce have
trivial limits. The purpose of this paper is to introduce for
certain classes of laminations $\mathcal{L}$ a
construction $\fgrm ({\cal L},x)$ called the fundamental germ, a
generalization of $\pi_{1}$ which represents an attempt to address this omission
in the theory of laminations.

The intuition which guides the construction is that of the
lamination as irrational manifold. Recall that for a pointed
manifold $(M,x)$, the deck group of the universal cover
$(\widetilde{M},\tilde{x})\rightarrow (M,x)$ -- which may be identified
with $\pi_{1}(M,x)$ -- reveals through its action how
to make identifications within $(\widetilde{M},\tilde{x})$ so as
to recover $(M,x)$ by quotient. Let us
imagine that we have disturbed the process of identifying
$\pi_{1}$ orbits, so that instead, points in an orbit merely
approximate one another through some auxiliary transversal space
$\textsf{T}$. The result is that $(\widetilde{M},\tilde{x})$ does
not produce a quotient manifold but rather coils upon itself, perhaps
forming a leaf $(L,x)$ of a lamination
$\mathcal{L}$. The germ of the transversal {\sf T} about $x$ may
be interpreted as the failed attempt of $(L,x)$ to form an
identification topology at $x$.  The fundamental germ $\fgrm({\cal
L},x)$ is then a device which records algebraically the dynamics
of $(L,x)$ as it approaches $x$ through the topology of $\sf T$.
See Figure 1.

\begin{figure}[htbp]
\centering \epsfig{file=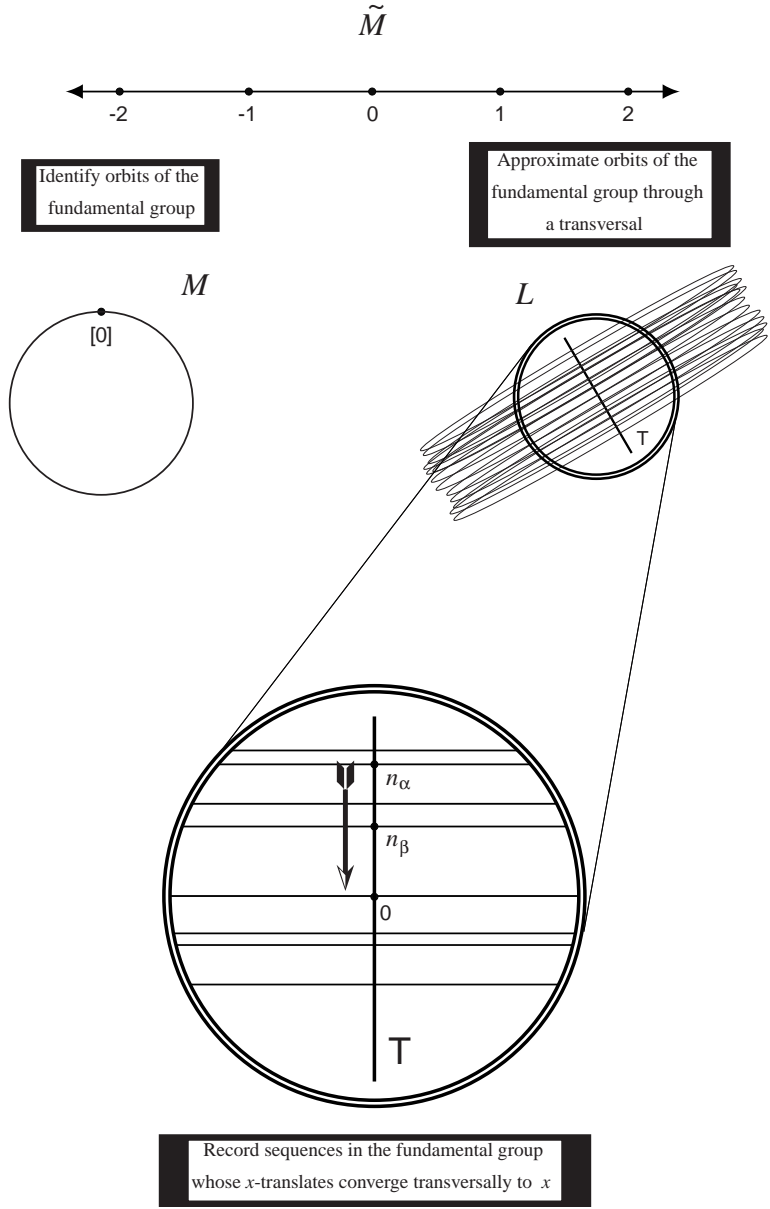, width=4in} \caption{The
Lamination as Irrational Manifold }
\end{figure}

One might define an element of $\fgrm ({\cal L},x)$ as a
tail equivalence class of a sequence of approaches $\{
x_{\alpha}\}$, where $L\ni x_{\alpha}\rightarrow x$ through $\sf
T$.  In this paper, the laminations under consideration (see \S 2) have the
property that there is a group $G$ acting on $L$ in such a way
that every approach is of the form $\{g_{\alpha}x\}$, for
$g_{\alpha}\in G$.  We then define $\fgrm ({\cal L},x)$ as the set
of tail equivalence classes of sequences of the
form $\{
g_{\alpha}h_{\alpha}^{-1}\}$, where $g_{\alpha}x,\;
h_{\alpha}x\rightarrow x$ in {\sf T}.  A groupoid structure on $\fgrm ({\cal L},x)$
is defined by component-wise multiplication of sequences, and
$\pi_{1}(L,x)$ is contained in $\fgrm ({\cal
L},x)$ as a subgroup.  In practice, $\fgrm ({\cal L},x)$ has no
additional structure; but for many reasonably well-behaved laminations such as
inverse limit solenoids, Sullivan solenoids and linear foliations
of torii, it is a group.  See \S\S 3 -- 7 for definitions
and examples.

When $\mathcal{L}=M$ is a manifold (a lamination with one leaf),
$\fgrm (M,x)$ is equal to ${}^{\ast}\pi_{1}(M,x)$, the
nonstandard version of $\pi_{1}(M,x)$: the group of tail
equivalence classes of {\em all} sequences in $\pi_{1}(M,x)$. When
$\mathcal{L}$ is a lamination contained in a manifold
$M$, under certain circumstances, \S 7, there is a map $\fgrm ({\cal L},x)\rightarrow
{}^{\ast}\pi_{1}(M,x)$ whose image consists of those classes of
sequences in $\pi_{1}(M,x)$ that correspond to the holonomy of
$\cal L$. Thus, in expanding $\pi_{1}$ to its nonstandard
counterpart, it is possible to detect -- algebraically -- sublaminations
invisible to $\pi_{1}$.

One can profitably think of $\fgrm ({\cal L},x)$ as made from
sequences of ``$G$-diophantine
approximations''. In the case of an irrational foliation
$\mathcal{F}_{r}$ of the torus $\T^{2}$ by lines of slope
$r\in\R\setminus\Q$, \S 4.4, this is literally true: the elements of
$\fgrm(\mathcal{F}_{r},x)$ are the equivalence classes of
diophantine approximations of $r$. More
generally, in $\fgrm$ one finds an algebraic-topological tool which
enables systematic translation of the geometry of laminations
into the algebra of (non-linear) diophantine approximation.

One can extend the definition
of the fundamental germ to include accumulations of $L$
on points of other leaves.  Thus if $\hat{x}$ is any point
of $\mathcal{L}$, we define $\fgrm (\mathcal{L},x,\hat{x})$
as the set of classes of sequences of the form $\{ g_{\alpha}\cdot h^{-1}_{\alpha}\}$
where $g_{\alpha}x,h_{\alpha}x\rightarrow\hat{x}$.
We suspect that, together with the
topological invariants of the leaves, the fundamental germs 
$\fgrm (\mathcal{L},x,\hat{x})$ will play a central role in the topological
classification of
laminations.

By unwrapping the accumulations of $L$ implied
by the fundamental germ $\fgrm (\mathcal{L},x)$, one obtains
the germ universal cover $\gcom{\mathcal{L}}$, \S 9, which is a kind
of nonstandard completion of $\widetilde{L}$.  If
$\fgrm (\mathcal{L},x)$ is a group, then under certain circumstances one 
may associate lamination coverings 
$\mathcal{L}_{\textsf{C}} :=\textsf{C}\backslash\gcom{\mathcal{L}}$ 
of $\mathcal{L}$ to every
conjugacy class of subgroup $\textsf{C}<\fgrm (\mathcal{L},x)$, and
when $\textsf{C}$ is a normal subgroup, the
quotient $\fgrm ({\cal L},x)/\textsf{C}$ may be identified with
the automorphism group of $\mathcal{L}_{\textsf{C}}\rightarrow
\mathcal{L}$. These considerations give rise to the beginnings of
a Galois theory of laminations, \S 10.

This first paper on the fundamental germ is foundational in nature.
One should not expect to find in it hard theorems, but rather
the description
of a complex and mysterious
object which reveals the explicit connection
between the geometry of laminations and the algebra of diophantine approximation.
Due to its somewhat elaborate construction, 
we shall confine ourselves here to the following
themes:
\begin{itemize}
\item Basic definitions: \S\S 1 -- 3.
\item Examples: \S\S 4 -- 7.
\item Functoriality: \S 8.
\item Covering space theory: \S\S 9,10.
\end{itemize}
The focus will be on
laminations which arise through group actions: suspensions,
quasi-suspensions, double coset foliations and locally-free Lie
group actions.  The exposition will be characterized by a
careful exploration of a number of concrete examples which serve
not only to illustrate the definitions in action but also to
indicate the richness of the algebra they produce.  In a second
installment \cite{Gendron2}, to appear elsewhere, the construction of $\fgrm$ 
will be extended to any lamination 
whose leaves admit a smooth structure.

\vspace{3mm}

\noindent {\bf Acknowledgements:}  I have benefited from
conversations with B. Le Roin, P. Makienko and especially A.
Verjovsky. In addition, the referee made many valuable suggestions
which helped to considerably improve the presentation.  
I would also like to thank the Instituto de
Matem\'{a}ticas (Cuernavaca) of the Universidad Nacional Aut\'{o}noma
de M\'{e}xico for providing generous financial support and a
pleasant work environment.

\section{Nonstandard Algebra}

We review facts concerning nonstandard
algebra, proofs of which may found in the literature. 
References: \ \cite{Goldblatt}, \cite{Robinson}.

Let $\N = \{ 0,1,2,\dots \}$, ${\mathfrak U}\subset 2^{\N}$ an ultrafilter
all of whose elements have
infinite cardinality. Given $\mathcal{S}=\{ S_{i}\}$ a sequence of
sets, write $ S_{X} = \prod_{j \in X} S_{j}.$
The {\em ultraproduct} is the direct limit
\[ [S_{i}] \;\; := \;\;
\lim_{\longrightarrow} S_{X}, \]
where the system maps are the cartesian projections.
If $S_{i}=S$ for all
$i$, the ultraproduct is called the {\em ultrapower} of $S$,
denoted $\bast S$.

If $\mathcal{S}$ consists of nested sets, denote by ${}^{\circledcirc}\mathcal{S}$
the set of sequences which converge with respect to $\mathcal{S}$. 
For each
$X\in\mathfrak{U}$, define a map $P_{X}:{}^{\circledcirc}\mathcal{S}
\rightarrow {}^{\circledcirc}\mathcal{S}$ by restriction of indices:
$ P_{X}\big(\{ x_{\alpha}\}\big)=\{ x_{\alpha}\}|_{\alpha\in X} $.
The {\em ultrascope} is the direct limit
\[ \bigodot S_{i} \;\; := \;\;
\lim_{\stackrel{\normalsize{\longrightarrow}}{P_{X}}}
{}^{\circledcirc}\mathcal{S} .\]
There is a
canonical inclusion
$ [S_{i}]\hookrightarrow
\bigodot S_{i}$ , and when $S_{i}=S$ for all $i$,
the ultrascope coincides with the
ultrapower.  
In general, we have
$\bigodot S_{i}=\bigcap \bast S_{i}\supseteq
\bast\left(\bigcap S_{i}\right)$,  where the
inclusion is an equality if and only if 
$S_{i}$ is eventually equal to a fixed set.

If $\mathcal{S}$ is a (nested) sequence of groups or rings, the
induced component-wise operations on sequences descend to
operations making the ultraproduct (the ultrascope) a
group or ring.  This is also true if $\mathcal{S}$ is a (nested)
sequence of fields: we remark here that the maximality
property of ultrafilters is required to rule out zero
divisors.  

If one uses a different ultrafilter $\mathfrak{U}'$ 
and if $\mathcal{S}$ is a (nested) sequence of groups, rings
or fields, then assuming the continuum hypothesis, it is classical
\cite{ChangKeisler} that the resulting ultraproduct is isomorphic to that formed
from $\mathfrak{U}$.  The same can shown for the ultrascope, however we
shall not pursue this point here. 

The ultrapower $\bast\R$ is called {\em nonstandard} $\R$. There
is a canonical embedding $\R\hookrightarrow\bast\R$ given by the
constant sequences, and we will not distinguish between $\R$ and
its image in $\bast \R$.  For $\bast x,\bast y\in \bast\R$, we write
$\bast x<\bast y$ if there exists $X\in\mathfrak{U}$ and
representative sequences $\{ x_{i}\}$, $\{ y_{i}\}$ such that
$x_{i} <y_{i}$ for all $i\in X$.  The non-negative nonstandard
reals are defined $\bast\R_{+}=\{ \bast x\in\bast\R
\, |\; \bast x\geq 0\}$. The Euclidean norm $|\cdot |$ on $\R$
extends to a $\bast\R_{+}$-valued norm on $\bast\R$. An
element $\bast x$ of $\bast\R$ is called {\em infinite} if for all
$r\in\R$, $|\bast x|>r$, otherwise $\bast x$ is called {\em finite}.  
$\bast\R$ is a totally-ordered, non-archimedian
field.

Here are two topologies that we may give $\bast\R$:
\begin{itemize}
\item The {\em enlargement topology} $\bast\tau$, generated by sets of the form 
$\bast A$,
where $A\subset\R$ is open. $\bast\tau$ is $2^{\rm
nd}$-countable but not Hausdorff.  
\item The {\em internal topology} $[\tau ]$, generated by sets of the form
$[A_{i}]$
where $A_{i}\subset \R$ is open for all $i$. $[\tau ]$ is Hausdorff 
but not $2^{\rm nd}$-countable.
\end{itemize}
We have $\bast\tau\subset [\tau]$, the inclusion being strict.
It is not difficult to see that $[\tau ]$ is just the order
topology.  

\begin{prop} $\big( \bast\R ,\; [\tau ]\big)$ is a real, infinite
dimensional
topological vector space.
\end{prop}

We note however that $\bast\R$ is not a topological group with
respect to $\bast\tau$.  Let $\bast\R_{\rm fin}$ be the set of 
finite elements of $\bast\R$.

\begin{prop}\label{subrng}
$\bast\R_{\rm fin}$ is a topological subring of $\bast\R$
with respect to both the $\bast\tau$ and $[\tau ]$ topologies.
\end{prop}

The set of {\em infinitesimals} is defined $\bast\R_{\epsilon}= \{ \bast\epsilon\,|\;
|\bast\epsilon |<M\text{ for all } M\in
\R_{+}\}$, a vector subspace of $\bast\R$.  
If $\bast x-\bast y\in\bast\R_{\epsilon}$,
we write $ \bast x\;\simeq\;\bast y$
and say that $\bast x$ is {\em infinitesimal} to $\bast y$.

\begin{prop}\label{quot}  $\bast\R_{\rm fin}$ is a local ring with maximal ideal
$\bast\R_{\epsilon}$ and \[ \bast\R_{\rm fin}/\bast\R_{\epsilon}
\;\;\cong\;\; \R ,\] a homeomorphism
with respect to the quotient $\bast\tau$-topology.
\end{prop}

We note that $\bast\R_{\epsilon}$ is clopen in the $[\tau ]$-topology; 
the quotient $[\tau ]$-topology
on $\bast\R_{\rm fin}/\bast\R_{\epsilon}$ is therefore discrete.
$\bast\R_{\epsilon}$ is not an ideal in $\bast\R$.
The vector space
\[ \bstar\R\;\; :=\;\; \bast\R/\bast\R_{\epsilon}, \]
equipped with the
quotient $\bast\tau$-topology, is called the {\em extended
reals}.  By Proposition~\ref{quot}, $\bstar\R$ contains a
subfield isomorphic to $\R$.

The results above
show that neither topology $\bast \tau$ or $[\tau ]$ can claim
to be preferred.  The lack of a canonical topology
on $\bast\R$ is a theme we will encounter again in \S 9, where
we will see that $\bstar\R$
may be viewed as the universal cover of a host of 1-dimensional
laminations, each one providing a different topology to $\bstar\R$
(and by pull-back to $\bast\R$).

Now let $\mathfrak{G}$ be any complete topological group.  Some of the
properties satisfied by $\bast\R$ also hold for
$\bast\mathfrak{G}$.  If $\tau$ denotes the topology of $\mathfrak{G}$, then the
topologies $\bast\tau$ and $[\tau ]$ are defined exactly as above.
$\bast\mathfrak{G}$ is a topological group in the
$[\tau ]$ topology, but not in the $\bast\tau$ topology.
Denote by $\bast\mathfrak{G}_{\epsilon}$ the classes of
sequences converging to the unit element 1.
$\bast\mathfrak{G}_{\epsilon}$ is a group since a product
of sequences converging to 1 in a topological group is again a
sequence converging to 1.
Let $\bast\mathfrak{G}_{\rm fin}$ be the subset of
$\bast\mathfrak{G}$ all of whose elements are represented by
sequences which converge to an element of $\mathfrak{G}$. We have
the following analogue of Proposition~\ref{quot}:

\begin{prop}\label{infinitesimalgroup}   $\bast\mathfrak{G}_{\epsilon}$ is a normal
subgroup of $\bast\mathfrak{G}_{\rm fin}$ and
\[ \bast\mathfrak{G}_{\rm fin}/\bast\mathfrak{G}_{\epsilon}
\;\;\cong\;\;\mathfrak{G},\] a homeomorphism with respect to the
quotient $\bast\tau$-topology.
\end{prop}

The left coset space
\[ \bstar\mathfrak{G}\;\; :=\;\;
\bast\mathfrak{G}/\bast\mathfrak{G}_{\epsilon}, \] 
with the
quotient $\bast\tau$-topology, is called the {\em extended}
$\mathfrak{G}$. It contains $\mathfrak{G}$ as a topological
subgroup.  If $\mathfrak{G}$ is compact or abelian, then
$\bstar\mathfrak{G}$ is a group, though in general it need not be.
We will avail ourselves of its natural structure as a 
$\bast \mathfrak{G}$-set with respect to the left muliplication action.

\section{Laminations Associated to Group Actions}

The laminations for which we shall define the fundamental
germ arise from actions of groups: we review them here
as a way of fixing notation. References:
\cite{CandelConlon}, \cite{Godbillon}, \cite{MooreSchochet}.

Let us begin by reviewing the definitions and terminology
surrounding the concept of a lamination.  
A {\em deck of cards} is a product $\R^{n}\times {\sf T}$, where
$\sf T$ is a topological space.  A {\em card} is a subset of the
form $C=O\times \{ \sf t\}$, where $O\subset \R^{n}$ is open and
${\sf t}\in {\sf T}$.  A {\em lamination} of dimension $n$ is a space $\mathcal{L}$ 
equipped
with a maximal atlas $\mathcal{A}=\{\phi_{\alpha}\}$ consisting of
charts with range in a fixed deck of cards $\R^{n}\times {\sf T}$,
such that each transition homeomorphism
$\phi_{\alpha\beta}=\phi_{\beta}\circ\phi_{\alpha}^{-1}$ satisfies
the following conditions:
\begin{enumerate}
\item For every card $C\in {\sf Dom}(\phi_{\alpha\beta})$,
$\;\phi_{\alpha\beta}(C)$ is a card.
\item The family of homeomorphisms $\{\phi_{\alpha\beta}(\cdot ,{\sf
t})\}$ is continuous in ${\sf t}$.\end{enumerate}
If {\sf T} is totally disconnected, we say that $\mathcal{L}$ is a {\em
solenoid}.

An open (closed) {\em transversal} in $\mathcal{L}$ is a subset of
the form $\phi_{\alpha}^{-1}(\{ x\}\times {\sf T}')$ where ${\sf
T}'$ is open (closed) in $\sf T$.  Note that an open (closed)
transversal need not be open (closed) in $\mathcal{L}$ {\em i.e.}
if $\mathcal{L}$ is a manifold (viewed as a trivial lamination) then every point
is an open transversal.  An open (closed) {\em flow box}
is a subset of the form $\phi_{\alpha}^{-1}(O\times {\sf T}')$,
where $O$ is open and ${\sf T}'\subset {\sf T}$ is open (closed).
A {\em plaque} in $\mathcal{L}$ is a subset of the form
$\phi_{\alpha}^{-1}(C)$ for $C$ a card in the deck $\R^{n}\times
{\sf T}$. A {\em leaf} $L\subset\mathcal{L}$ is a maximal
continuation of overlapping plaques in $\mathcal{L}$. Note that $\mathcal{L}$
is the disjoint union of its leaves; we denote by $L_{x}$ the leaf
containing the point $x$.
A lamination
is {\em weakly minimal} if it has a dense leaf; it is {\em
minimal} if all of its leaves are dense.   
A transversal which meets every leaf is called {\em complete}.
Unless we say otherwise, all transversals in this paper will 
be complete and open.  Two laminations $\mathcal{L}$ and $\mathcal{L}'$
are said to be {\em homeomorphic} if there  
is a homeomorphism $f:\mathcal{L}\rightarrow\mathcal{L}'$ mapping
leaves homeomorphically onto leaves and transversals homeomorphically onto 
transversals.

\subsection{Suspensions}\label{suspensions}

Let $B$ be a manifold in which $\pi_{1}B$ acts without fixed points.
Let $F$ be a topological
space and $\rho :\pi_{1}B\rightarrow {\sf Homeo}(F)$ a
representation. The {\em suspension} of $\rho$ is the space
\[ \mathcal{L}_{\rho}\;\; =\;\; \widetilde{B}\times_{\rho} F \]
defined by quotienting $\widetilde{B}\times F$ by the
diagonal action of $\pi_{1}B\,$, $\alpha\cdot
(\tilde{x},\, t) = 
(\alpha\cdot\tilde{x},\,\rho_{\alpha }(t))$.  The suspension is
a fiber bundle over $B$ with model fiber $F$.
If $F=\mathfrak{G}$ is a topological group and 
$\varphi:\pi_{1}B\rightarrow \mathfrak{G}$ a homomorphism, then the
representation $\rho :\pi_{1}B\rightarrow {\sf Homeo}(\mathfrak{G})$
defined $\rho_{\gamma}(g)=g\cdot \varphi(\gamma^{-1})$ gives rise
to what we call a $\mathfrak{G}$-{\em suspension}, denoted $\mathcal{L}_{\varphi}$, 
a principle $\mathfrak{G}$-bundle over $B$. 

The action of $\pi_{1}B$ used to define $\mathcal{L}_{\rho}$ 
is properly discontinuous and leaf preserving, hence $\mathcal{L}_{\rho}$
is a lamination modeled on the deck of cards $\widetilde{B}\times F$.
If $K=\ker (\rho )$ and $(L,x)\subset\mathcal{L}_{\rho}$ is a pointed leaf, we have
$K\unlhd\pi_{1}(L,x)$.  $\mathcal{L}_{\rho}$ is minimal (weakly-minimal) if and only
if every (at least one) $\rho (\pi_{1}B )$ orbit is dense.

 The restriction $p|_{L}$ of the projection $p:\mathcal{L}_{\rho}\rightarrow
B$ to a leaf $L$ is a covering map.  
Suppose that $p_{L}$ is a Galois covering (we say
that $L$ is {\em Galois}).  The deck group $D_{L}$ of $p|_{L}$ has
the property that
\[D_{L}\cdot x \;\;=\;\; L\cap F_{x},\] where $F_{x}$ is the fiber of $p$
through $x$. In particular, if we give $(L\cap F_{x})\subset
F_{x}$ the subspace topology, we have an inclusion
\[ D_{L}\hookrightarrow {\sf Homeo}(L\cap F_{x}).\]

A manifold $B$ is a suspension with $F$ a point
and
$\rho:\pi_{1}B\rightarrow F$
trivial.
The following subsections discuss examples which are more interesting.

\subsubsection{Inverse Limit Solenoids}

Let $\mathcal{C}=\big\{  \rho_{\alpha}:M_{\alpha}\rightarrow M
\big\}$ be an inverse system of pointed manifolds and Galois
covering maps with initial object $M$; denote by
\[ \widehat{M}\;\;=\;\;\widehat{M}_{\mathcal{C}}\;\; :=\;\;
\lim_{\longleftarrow}M_{\alpha}\] the limit.  By definition
$\widehat{M}\subset\prod M_{\alpha}\,$, so elements of
$\widehat{M}$ are denoted $\hat{x}=(x_{\alpha})$, where
$x_{\alpha}\in M_{\alpha}$.  The natural projection onto the base
surface is denoted $p :\widehat{M}\rightarrow M$.
We may identify the universal covers
$\widetilde{M}_{\alpha}$ with $\widetilde{M}$ and choose the
universal covering maps $\widetilde{M}\rightarrow M_{\alpha}$ to
be compatible with the system $\mathcal{C}$.  By universality,
there exists a canonical map
$i :\widetilde{M}\longrightarrow \widehat{M}$.

Let $H_{\alpha}=(\rho_{\alpha})_{\ast}(\pi_{1}M_{\alpha})<\pi_{1}M$.
Associated to $\mathcal{C}$ is the inverse limit of deck groups
\[ \hat{\pi}_{1}M \;\; :=\;\; \lim_{\longleftarrow} \pi_{1}M/H_{\alpha} ,
\]
a Cantor group since the $\pi_{1}M/H_{\alpha}$ are finite.
By universality of inverse limits, the projections
$\pi_{1}M\rightarrow\pi_{1}M/H_{\alpha}$ yield a canonical
homomorphism
$\iota :\pi_{1}M\longrightarrow \hat{\pi}_{1}M$
with dense image.   The
closures of the images $\iota (H_{\alpha})$ are clopen, and give a
neighborhood basis about 1.  
Let $\mathcal{L}_{\iota}$ be the associated $\hat{\pi}_{1}M$-suspension.

\begin{prop}  $\widehat{M}$ is homeomorphic to $\mathcal{L}_{\iota}$.
In particular, $\widehat{M}$ is a solenoid.
\end{prop}

\begin{proof} Let $\Upsilon: \widetilde{M}\times\hat{\pi}_{1}M\rightarrow
\widehat{M}$ be the map defined $(\tilde{x},\hat{g})\mapsto
\hat{g}\cdot i(\tilde{x})$. $\Upsilon$ is invariant with respect
to the diagonal action of $\pi_{1}M$, and descends to a
homeomorphism $
\widetilde{M}\times_{\rho}\hat{\pi}_{1}M\rightarrow\widehat{M}$.
\end{proof}

\subsubsection{Linear Foliations of Torii}

Let $V$ be a $p$-dimensional subspace of $\R^{p+q}$. Denote by
$\widetilde{\mathcal{F}}_{V}$ the foliation of $\R^{p+q}$ by cosets
$\textbf{v}+V$.  The image $\mathcal{F}_{V}$ of $\widetilde{\mathcal{F}}_{V}$ 
in the torus $\T^{p+q}=\R^{p+q}/\Z^{p+q}$ gives a
foliation of the latter by Euclidean manifolds.
$V$ may be regarded as the graph of a $q\times p$ matrix map
\[\mathbf{R}: \R^{p}\rightarrow\R^{q}\] whose columns are 
independent.  For $\textbf{y}\in \R^{q}$, denote by
$\overline{\textbf{y}}$ its image in $\T^{q}$.
Let $\varphi_{\textbf{R}}:\Z^{p}\rightarrow \T^{q}$ be the
homomorphism defined 
\[ \varphi_{\textbf{R}}(\textbf{n})\;\; =\;\;
\overline{\textbf{R}\textbf{n}}, \]
and denote by $\mathcal{L}_{\varphi_{\textbf{R}}}$ the corresponding
$\T^{q}$-suspension.

\begin{prop}\label{toral} 
$\mathcal{F}_{V}$ is homeomorphic to {\rm $\mathcal{L}_{\varphi_{\textbf{R}}}$}.
\end{prop}

\begin{proof} Let $P_{0}:\R^{p+q}=\R^{p}\times\R^{q}\rightarrow \R^{p}\times \T^{q}$ be the
map defined $(\textbf{x},\, \textbf{y})
\mapsto \left(\textbf{x},\,
\overline{\textbf{y}}-\overline{\textbf{R}\textbf{x}}\right)$.  
Let $P$ be the composition of
$P_{0}$ with the projection $\xi :\R^{p}\times
\T^{q}\rightarrow \mathcal{L}_{\varphi_{\textbf{R}}}$. Then $P$ is a
covering homomorphism with kernel $\Z^{p+q}$, hence
$\mathcal{L}_{\varphi_{\textbf{R}}}\approx\T^{p+q}$. Since
$V=(\textbf{x},\textbf{R}\textbf{x})$, we have $P(V) =\xi
(\R^{p}\times \overline{\textbf{0}});$ thus $P(V)$ is a leaf
of the suspension. It follows that $P$ defines a map
$\widetilde{\mathcal{F}}_{V}\rightarrow \mathcal{L}_{\varphi_{\textbf{R}}}$
which descends to the desired homeomorphism.
\end{proof}

Let $\textbf{r}_{i}$ be the $i$th column vector of $\mathbf{R}$.
If $\textbf{r}_{i}\in\Q^{q}$ for all $i$, the leaves of
$\mathcal{F}_{V}$ are homeomorphic to $\T^{p}$ and are not dense.
If at least one of the $\textbf{r}_{i}$ has an irrational
coordinate, then the leaves of ${\cal F}_{V}$ are non-compact and
dense, homeomorphic to the quotient of $\R^{p}$ by a discrete subgroup with
as many generators as rational $\textbf{r}_{i}$.

\subsubsection{Anosov Foliations}\label{Anosov}

Let $\Sigma=\HP^{2}/\Gamma$ be a hyperbolic surface and
let $\rho :\Gamma\rightarrow {\sf Homeo} (\SC^{1})$ be defined by
extending the action of $\Gamma$ on $\HP^{2}$ to
$\partial\HP^{2}\approx \SC^{1}$. The suspension 
\[\mathcal{F}_{\Gamma}=\HP^{2}\times_{\rho} \SC^{1}\] is called an {\em Anosov
foliation}.  Note that $\mathcal{F}_{\Gamma}$ is not an $\SC^{1}$-suspension.
It is classical that the underlying space of $\mathcal{F}_{\Gamma}$ 
is homeomorphic to the unit tangent bundle ${\rm
T}^{1}_{\ast}\Sigma$.

\subsection{Quasisuspensions}\label{quasisuspensions}

Let $\mathcal{L}_{\rho}=\widetilde{B}\times_{\rho}F$ be a suspension over
a manifold $B$. We say that $\mathcal{L}_{\rho}$ is {\em Galois} if every
leaf of $\mathcal{L}_{\rho}$ is Galois.
Throughout this section, $\mathcal{L}_{\rho}$ will be a Galois
suspension.  We define an action of $\pi_{1}B$ on $\mathcal{L}_{\rho}$ by
\[ x\;\longmapsto\; \bar{\gamma}\cdot x,\] where, for $x$
contained in the leaf $L$, $\bar{\gamma}$ is the image of
$\gamma\in\pi_{1}B$ in $\pi_{1}B/(p_{L})_{\ast}(\pi_{1}L)$ $\cong$
$D_{L}$ = the deck group of $p|_{L}$.  

Let $\mathcal{X}\subset \mathcal{L}_{\rho}$ be any closed subset which
is invariant with respect to the action of $\pi_{1}B$. 
Let $\mathcal{L}_{0}:=\mathcal{L}_{\rho}\setminus \mathcal{X}$, which 
is a lamination mapping to $B$.  If $\mathcal{L}_{\rho}$ is minimal, then 
$\mathcal{X}$ is the preimage of a subset $X\subset B$, hence $\mathcal{L}_{0}$
is a fiber bundle over $B_{0}=B\setminus X$.  In general, we shall
define the fibers of $\mathcal{L}_{0}$ over $x\in B$ to be the preimages 
of the map $\mathcal{L}_{0}\rightarrow B$.

A lamination homeomorphism $f:\mathcal{L}_{0}\rightarrow \mathcal{L}_{0}$ is
{\em weakly fiber-preserving} if for every fiber $F_{x}$ over $B$,
\begin{equation}\label{wkfib} f(F_{x})\;\;=\;\;\bigcup_{i=1}^{n}
E_{x_{i}},\end{equation} where $E_{x}\subset F_{x}$ denotes a
subset of the fiber $F_{x}$. The collection ${\sf
Homeo}_{\omega{\sf -fib}}(\mathcal{L}_{0})$ of weakly fiber-preserving
homeomorphisms is clearly a group.
Since the fibers are disjoint, each $E_{x_{i}}$ occurring in
(\ref{wkfib}) must be open in $F_{x_{i}}$.  In particular, if the fibers
are connected, a weakly fiber-preserving homeomorphism is
fiber-preserving. Thus, the concept of a weakly fiber-preserving
homeomorphism differs from that of a fiber-preserving homeomorphism
when the fibers are 
disconnected {\em e.g.} when $\mathcal{L}_{0}$ is a solenoid.

\begin{defi} Let $\mathcal{L}_{0}$ be as above and suppose 
$H<{\sf Homeo}_{\omega{\sf -fib}}(\mathcal{L}_{0})$ is a subgroup
acting properly discontinuously on $\mathcal{L}_{0}$. The quotient 
\[\mathcal{Q} \;\;=\;\; H\backslash \mathcal{L}_{0}\] is a lamination called a {\bf
quasisuspension} (over $B$).  
\end{defi}

We consider now two examples.

\subsubsection{The Sullivan Solenoid}  
The following important example comes from holomorphic dynamics.
Let $U,V\subset \C$ be regions conformal to the unit disc, with
$\overline{U}\subset V$. Recall that a {\em polynomial-like map} is
a proper conformal map $f:U\rightarrow V$.  The conjugacy class of $f$ is
uniquely determined by a pair $(p,\partial f)$,
where $p$ is a complex polynomial of degree $d$ and $\partial f:
\SC^{1}\rightarrow \SC^{1}$ is a smooth, expanding map of degree $d$ 
\cite{DouadyHubbard}.  The space
\begin{equation}\label{defininginvsys}\SO\;\;=\;\;\lim_{\longleftarrow}
\Big( \SC^{1}\stackrel{\partial
f}{\longleftarrow
}\SC^{1}\stackrel{\partial f}{\longleftarrow}
\SC^{1}\stackrel{\partial f}{\longleftarrow}\cdots \Big)
\end{equation} is an inverse
limit solenoid which may be identified with the $\widehat{\Z}_{d}$-suspension
$\mathcal{L}_{\iota}=\R\times_{\rho}\widehat{\Z}_{d}$,
where $\widehat{\Z}_{d}$ is the group of $d$-adic integers and
$\imath :\Z\hookrightarrow \widehat{\Z}_{d}$ is the canonical inclusion. Every
leaf of $\SO$ is homeomorphic to $\R$. $\partial f$
defines a self map of the inverse system in (\ref{defininginvsys}), inducing a
homeomorphism $\partial \hat{f}:\SO
\rightarrow\SO$.

Consider the suspension
\[ \widehat{\D}\;\; =\;\;\HP^{2}\times_{\rho}\widehat{\Z}_{d}\]
obtained by extending to $\HP^{2}\times\widehat{\Z}_{d}$ the 
identification used to define $\mathcal{L}_{\iota}$
{\em e.g.} $(z,\hat{n})\sim (\gamma^{m}\cdot
z,\;\rho_{m}(\hat{n}))$ for $m\in\Z$, where
\[ \gamma \;\; =\;\;\left(
\begin{array}{cc}
           1 &  1 \\
            0 &  1
          \end{array}
                       \right)\]
is the affine extension of the map $x\mapsto x+1$ to $\HP^{2}$.
The base of the suspension $\widehat{\D}$
is the punctured hyperbolic disc $\D^{\ast}=\langle \gamma \rangle
\backslash\HP^{2}$, and its ideal boundary may be identified
with $\SO$.

The map $\partial \hat{f}$ extends to a weakly fiber-preserving
homeomorphism $\hat{f}:\widehat{\D}\rightarrow\widehat{\D}$ which
acts properly discontinuously on $\widehat{\D}$.  The quotient
\[ \widehat{\D}_{f} \;\; :=\;\; \langle \hat{f}\rangle\backslash
\widehat{\D} \] is a quasisuspension called the {\em Sullivan
solenoid} \cite{Sullivan}.

\subsubsection{The Reeb Foliation}\label{ReebFoliation}  Let $\R_{+}=[0,\infty )$,
consider the
trivial suspension $\C\times \R_{+}$ over $\C$, and 
denote $(\C\times \R_{+})^{\ast}=\C\times \R_{+}\setminus \{ (0,0)\}$.  
Fix $(\mu ,\lambda )\in (\C\times \R_{+})^{\ast}$ with 
$| \mu |,\lambda >1$, $\mu\not=\lambda$.
Then multiplication by $(\mu ,\lambda )$ in $(\C\times \R_{+})^{\ast}$ is
a fiber-preserving lamination homeomorphism giving rise to an action by $\Z$.
The resulting quasisuspension
\[ \mathcal{F}_{\sf Reeb} \;\;=\;\; \Z\backslash (\C\times \R_{+})^{\ast}\]
has underlying space a solid torus, and is called the {\em Reeb foliation}.

Let $P:(\C\times\R_{+})^{\ast}\rightarrow\mathcal{F}_{\sf Reeb}$
denote the projection map. The leaves of $\mathcal{F}_{\sf Reeb}$
are of the form:
\begin{enumerate}
\item $L_{t}=P(\C\times\{ t\})\cong\C$, for $t>0$.
\item $L_{0}=P(\C^{\ast}\times\{ 0\})
\cong \C^{\ast}/<\mu >$.  
\end{enumerate}
The fiber tranversals of $\mathcal{F}_{\sf Reeb}$ are of the form:
\begin{enumerate}
\item $T_{z}= P(\{ z\}\times \R_{+})\approx\R_{+}$, $z>0$.  Every leaf of
$\mathcal{F}_{\sf Reeb}$ intersects $T_{z}$.
\item $T_{0}=P(\{ 0\}\times (0,\infty ))\approx \SC^{1}$.  Every leaf except $L_{0}$
intersects $T_{0}$.
\end{enumerate}

There is an action of $\Z$ on $\mathcal{F}_{\sf Reeb}$ 
induced by the map $(z,\; t)\mapsto (\mu^{n}z,\; t)$.  
For $x\in \mathcal{F}_{\sf Reeb}$, we write this action $x\mapsto n\cdot x$.  
For every $t$
we have $n\cdot L_{t}$ = $L_{t}$ and for all $z$,
$n\cdot T_{z}$ = $T_{z}$.  Note that this action is the identity on $L_{0}$.

\subsection{Double Coset Foliations}\label{DoubleCoset}

Let $\mathfrak{G}$ be a Lie group, $\mathfrak{H}$ a closed Lie
subgroup, $\Gamma <\mathfrak{G}$ a discrete subgroup.  The
foliation of $\mathfrak{G}$ by right cosets $\mathfrak{H}g$
descends to a foliation $\mathcal{F}_{\mathfrak{H},\Gamma}$ of
$\mathfrak{G}/\Gamma$, called a {\em double coset foliation}.

Let $\Gamma$ be a co-finite volume Fuchsian group.  Denote by
$\Sigma =\HP^{2}/\Gamma$ and by 
${\rm T}^{1}_{\ast}\Sigma$ the unit tangent bundle of $\Sigma$. Recall
that every $v\in {\rm T}_{\ast}^{1}\HP^{2}$ determines three
oriented, parametrized curves: a geodesic $\gamma$ and two
horocycles $\mathfrak{h}_{+}$, $\mathfrak{h}_{-}$ tangent to,
respectively, $\gamma\, (\infty )$ and $\gamma\, (-\infty )$. By
parallel translating $v$ along these curves, we obtain three flows
on ${\rm T}_{\ast}^{1}\HP^{2}$. The three flows are
$\Gamma$-invariant, and define flows on ${\rm
T}^{1}_{\ast}\Sigma$.  The corresponding foliations are denoted
${\sf Geod}_{\Gamma}$, ${\sf Hor}^{+}_{\Gamma}$ and ${\sf
Hor}^{-}_{\Gamma}$.

Now let $\mathfrak{G}=SL(2,\R )$ and take $\mathfrak{H}$ to be one
of the 1-parameter subgroups $H^{+}=\{ A^{+}_{r}\}$, $H^{-}=\{
A^{-}_{r}\}$ and $G=\{ B_{r}\}$, where
\[ A^{+}_{r}=\left(
\begin{array}{cc}
1 & r \\
0 & 1
\end{array}
\right),\;\;\;\;  A^{-}_{r}=\left(
\begin{array}{cc}
1 & 0 \\
r & 1
\end{array}
\right)\;\;\;\text{and}\;\;\; B_{r}=\left(
\begin{array}{cc}
e^{r/2} & 0 \\
0 & e^{-r/2}
\end{array}
\right)    \] for $r\in\R$.
Then it is classical that the foliations $\mathcal{F}_{G,\Gamma}$ and
$\mathcal{F}_{H^{\pm},\Gamma}$ are homeomorphic to ${\sf
Geod}_{\Gamma}$ and ${\sf Hor}^{\pm}_{\Gamma}$, respectively.
Note also that the Anosov foliation $\mathcal{F}_{\Gamma}$
is homeomorphic to the sum ${\sf Geod}_{\Gamma}\oplus {\sf
Hor}^{+}_{\Gamma}$.

\subsection{Locally-Free Lie Group Actions}\label{locallyfree}

Let $\mathfrak{B}$ be a Lie group of dimension $k$, $M^{n}$ an
$n$-manifold, $n>k$, $X$ a subspace of $M^{n}$.  A continuous
representation $\theta :\mathfrak{B}\rightarrow {\sf Homeo}(X)$ is
called {\em locally free} if for all $x\in X$, the isotropy
subgroup $I_{x}<\mathfrak{B}$ is discrete.  If for any pair $x,y\in X$,
their $\mathfrak{B}$-orbits are either disjoint or coincide, 
then $X$ has the structure
of a lamination $\mathcal{L}_{\mathfrak{B}}$ whose leaves are the
$\mathfrak{B}$-orbits.

For example, let $M^{n}$ be a Riemannian manifold.  Fix a tangent
vector $v\in {\rm T}_{x}M$.  Let $l\subset M^{n}$ be the
complete geodesic determined by $v$, $X$ its closure (itself a
union of geodesics). Then there is a locally free action of $\R$
given by geodesic flow along $X$, and $X$ is a lamination when $l$ is simple. 
When $M^{n}=\Sigma$ is a hyperbolic
surface and $l$ is simple,
we obtain a geodesic lamination in $\Sigma$ in
the sense of \cite{Thurston}, a solenoid since its transversals are 
totally-disconnected.

\section{The Fundamental Germ}  

Let $\mathcal{L}$ be any of the laminations considered
in the previous section and let $L\subset\mathcal{L}$
be a fixed leaf. If $\mathcal{L}=H\backslash\mathcal{L}_{0}$ is a quasisuspension
let $L_{0}\subset \mathcal{L}_{0}$ be a leaf lying over $L$.
The {\em diophantine group} $G_{L}$ of $\mathcal{L}$
with respect to $L$ is
\begin{itemize}
\item  $\pi_{1}B$ if $\mathcal{L}$ is a suspension.
\item  The group generated by $\pi_{1}B$, $H_{L}=\{ h\in H\, |\; h(L_{0})=L_{0}\}$ and 
$\pi_{1}L$
(viewed as groups acting on $\widetilde{L}$) if $\mathcal{L}$ is a quasisuspension.
\item  The group $\widetilde{\mathfrak{H}}$ if $\mathcal{L}$ is a double
coset.
\item  The group $\widetilde{\mathfrak{B}}$ if $\mathcal{L}$ is a locally
free Lie group action.
\end{itemize}
Note that in every case, $\pi_{1}L< G_{L}$.

Let $\hat{x}\in\mathcal{L}$ and $T$ a transversal
containing $\hat{x}$.  Denote by $\widetilde{T}_{L}\subset\widetilde{L}$  
the set of points lying over
$T\cap L$.  Then $T$ is said to be a {\em diophantine transversal} if for every
leaf $L$ and $\tilde{x}\in\widetilde{T}_{L}$, 
any $\tilde{y}\in\widetilde{T}_{L}$ may be written in the
form $\tilde{y}=g\cdot\tilde{x}$ for some $g\in G_{L}$.  For $\tilde{x}\in\widetilde{T}_{L}$
fixed, we call $\{ g_{\alpha}\}\subset G_{L}$ a {\em $G_{L}$-diophantine
approximation} of $\hat{x}$ along $T$ based at $\tilde{x}$ 
if $\{ g_{\alpha}\cdot\tilde{x}\}$
projects in $L$ to a sequence converging to $\hat{x}$ in $T$. 
The image of all such
$G_{L}$-diophantine approximations in $\bast G_{L}$ is denoted 
\[ \bast {\sf D} (\tilde{x}, \hat{x}, T),\]
and when $\hat{x}=x$ we write $\bast {\sf D} (\tilde{x}, T)$.
If there
are no $G_{L}$-diophantine approximations of $\hat{x}$ along $T$ based
at $\tilde{x}$, we define $\bast {\sf D} (\tilde{x}, \hat{x}, T)=0$.
Note that if $\tilde{x}'=\gamma\cdot\tilde{x}$ for $\gamma\in\pi_{1}L <G_{L}$ 
then 
\begin{equation}\label{changeoflift}
\bast {\sf D} (\tilde{x}', \hat{x}, T)\cdot\gamma \;\; =\;\;
\bast {\sf D} (\tilde{x}, \hat{x}, T).
\end{equation}
Let 
$\bast {\sf D}(\tilde{x},\hat{x},T)^{-1}$ consist of the set of inverses $\bast g^{-1}$
of classes belonging to $\bast {\sf D}(\tilde{x},\hat{x}, T)$.

\begin{defi}\label{FGDefn}  Let $\mathcal{L}$, $L$, $x$, $\hat{x}$ and $T$ be as above. 
The {\bf fundamental
germ} of $\mathcal{L}$ based at $\hat{x}$ along $x$ and $T$ is
\[ \fgrm (\mathcal{L},x, \hat{x}, T)\;\; =\;\; \bast {\sf D}(\tilde{x},\hat{x},T)\cdot
\bast {\sf D}(\tilde{x},\hat{x},T)^{-1} \]
where $\tilde{x}$ is any point in $\widetilde{L}$ lying over $x$.
\end{defi}

By (\ref{changeoflift}),  
$\fgrm (\mathcal{L},x, \hat{x}, T)$ does not depend on the choice of
$\tilde{x}$ over $x$.  When $x=\hat{x}\in L$,
we write $\fgrm (\mathcal{L},x, T)$.  Observe in this case that
$\gzg{\pi}_{1}(\mathcal{L},x,T)$ contains a
subgroup isomorphic to $\bast\pi_{1}(L,x)$. 

We now describe a groupoid structure on $\fgrm (\mathcal{L},x, \hat{x}, T)$ .
To do this, we define a unit space on which it acts: let $\bstar 
{\sf D} (\tilde{x}, \hat{x}, T)$ be the image of 
$\bast {\sf D} (\tilde{x}, \hat{x}, T)$
in $\bstar G_{L}$, for any $\tilde{x}$ over $x$.
We say that $\bast u\in \fgrm (\mathcal{L},x, \hat{x}, T)$ is
defined on $\bstar g\in \bstar {\sf D} (\tilde{x}, \hat{x}, T)$ if 
$\bast u\cdot \bstar g\in \bstar {\sf D} (\tilde{x}, \hat{x}, T)$. 
Here we are using the left action of $\bast G_{L}$ on $\bstar G_{L}$. 
Having defined the domain and range of elements of 
$ \fgrm (\mathcal{L},x, \hat{x}, T)$,
it is easy to see that 
$\fgrm (\mathcal{L},x, \hat{x}, T)$ is a groupoid,
as every element 
has an inverse by construction.  This groupoid structure
does not depend on the choice of $\tilde{x}$ over $x$.

\section{The Fundamental Germ of a Suspension}

In the case of a suspension $\mathcal{L}_{\rho}=\widetilde{B}\times_{\rho}F$, 
any fiber over the base $B$
is a diophantine transversal.  Conversely, any 
diophantine transversal is an open subset
of a fiber transversal.  It follows that any
two diophantine transversals $T, T'$ through a given point $\hat{x}$
define the same set of $G_{L}$-diophantine approximations.  
Thus

\begin{prop}\label{nodependenceonT}  If $T$ and $T'$ are diophantine 
transversals through $\hat{x}$ then
\[\fgrm (\mathcal{L}_{\rho}, x, \hat{x}, T)\;\; =\;\;
\fgrm (\mathcal{L}_{\rho}, x, \hat{x}, T').\]
\end{prop}

Accordingly for suspensions we drop mention of the transversal and write 
$\fgrm (\mathcal{L},x, \hat{x})$.  
We note that since the diophantine group $G_{L}=\pi_{1}B$ is discrete, 
$\bast G_{L}=\bstar G_{L}$
and the unit space for the groupoid structure is just
$\bast {\sf D}(\tilde{x}, \hat{x})$.

\subsection{Manifolds}

A manifold is a lamination with just one leaf, which 
can be viewed as the suspension
of the trivial representation of its fundamental group.  Since
a fiber transversal is just a point, we have immediately   

\begin{prop}  If $M$ is a manifold then
\[ \fgrm (M,x)\;\; =\;\; \bast {\sf D}(\tilde{x})\;\; =\;\;\bast \pi_{1}(M,x).\]
\end{prop}

\subsection{$\mathfrak{G}$-Suspensions}

Let $\varphi :\pi_{1}B\rightarrow \mathfrak{G}$ be a homomorphism, 
$\mathcal{L}_{\varphi}$ the corresponding $\mathfrak{G}$-suspension.
Let $\{U_{i}\}$ be a neighborhood basis about 1 in $\mathfrak{G}$ and
define a collection of nested sets $\{ G_{i}\}$ by 
$G_{i}=\{ \gamma\in\pi_{1}B\, |\; h(\gamma )\in U_{i}\}$.  Note
that the ultrascope $\bigodot G_{i}$ is a subgroup of $\bast\pi_{1}B$.
In fact, if $\bast\varphi :\bast\pi_{1}B\rightarrow\bast \mathfrak{G}$
is the nonstandard version of $\bast\varphi$, then 
\[ \bigodot G_{i}\;\;=\;\;\bast\varphi^{-1}(\bast\mathfrak{G}_{\epsilon}).\]

\begin{theo}\label{fgrmGsusp}  If $\varphi$ has dense image, then for any pair
$x,\hat{x}$ belonging to a diophantine transversal,
$\fgrm (\mathcal{L}_{\varphi},x,\hat{x} )$ is a group isomorphic to
$\bigodot G_{i}$.  
\end{theo}

\begin{proof} Let $\bast g\in \bast {\sf D}(\tilde{x}, \hat{x})$.
Then any other element $\bast g'\in \bast {\sf D}(\tilde{x}, \hat{x})$ 
may be written in the form $\bast g\cdot\bast h$ where $\bast h\in \bigodot G_{i}$.
It follows immediately that 
\[ \fgrm (\mathcal{L}_{\varphi},x,\hat{x})
\;\; =\;\; \bast g \cdot \big(\bigodot G_{i}\big) \cdot \bast g^{-1}\;\;
\cong\;\;\bigodot G_{i}.\]
Because the unit space $\bast {\sf D}(\tilde{x}, \hat{x})$ is invariant
under left-multiplication by its elements, it follows that
$\fgrm (\mathcal{L}_{\varphi},x,\hat{x})$ acts on it as a group, its groupoid law
coinciding with multiplication in $\bigodot G_{i}$.
\end{proof}

For $\mathfrak{G}$-suspensions with $\varphi$ having dense image, 
we can thus reduce our notation to 
$\fgrm (\mathcal{L}_{\varphi})$.

Let $\bast\varphi:\bast \pi_{1}B\rightarrow\bast \mathfrak{G}$
be the induced map of nonstandard groups, and denote by
$\bast\pi_{1}B_{\rm fin}$ the subgroup 
$\bast\varphi^{-1}(\bast\mathfrak{G}_{\rm fin})$.  The following theorem
can be used to display many familiar topological groups as algebraic quotients
of nonstandard versions of discrete groups.  

\begin{theo}\label{nonstdquot}  If $\varphi$ has dense image, 
then $\fgrm (\mathcal{L}_{\varphi})$ is a normal subgroup of 
$\bast\pi_{1}B_{\rm fin}$ with
\[ \bast\pi_{1}B_{\rm fin}\big/\fgrm (\mathcal{L}_{\varphi})\;\;\cong\;\; \mathfrak{G}.\]
\end{theo}

\begin{proof} Since $\varphi$ has dense image, the composition of homomorphisms
$\bast\pi_{1}B_{\rm fin}\rightarrow \bast\mathfrak{G}_{\rm fin}\rightarrow
\mathfrak{G}$ -- where the first arrow is $\bast\varphi$ -- is surjective with kernel 
$\bast\varphi^{-1}(\bast\mathfrak{G}_{\epsilon})=\fgrm (\mathcal{L}_{\varphi})$.
\end{proof}

\subsection{Inverse Limit Solenoids}\label{invlimit}

Let $\widehat{M}$ be an inverse limit solenoid over the base $M$,
and let $\{ H_{i}\}$ be a sequence of subgroups of $\pi_{1}M$
cofinal in the collection of subgroups in the defining inverse system.  
By the discussion
in \S 2.1.1, the collection of closures $\{ \widehat{H}_{i}\}
\subset\hat{\pi}_{1}M$ defines a neighborhood basis about $1$. 
Since $\widehat{M}$ is a $\hat{\pi}_{1}M$-suspension in which
$\varphi$ is dense, it follows from Theorem~\ref{fgrmGsusp} that
$\fgrm (\widehat{M},x,\hat{x})$ is a group
isomorphic to $\bigodot H_{i}$.

For example, consider a solenoid $\widehat{\SC}$
over $\SC^{1}$.  Here, each $H_{i}$ is an ideal in $\Z$, hence
$\fgrm \widehat{\SC}$ is an ideal in the ring $\bast\Z$ =
nonstandard $\Z$.   When $H_{i}=(d^{i})$ for $d\in\Z$ fixed, we denote the resulting
germ $\bast\Z_{\hat{\epsilon}}(d)$ and when $H_{i}= (i)$
we write $\bast\Z_{\hat{\epsilon}}$.  Being uncountable, 
these ideals are
not principal, so $\bast\Z$, unlike $\Z$, in not a PID.
By Theorem~\ref{nonstdquot}, we have $\bast\Z/\bast\Z_{\hat{\epsilon}}\cong
\widehat{\Z}$ and $\bast\Z/\bast\Z_{\hat{\epsilon}}(d)\cong
\widehat{\Z}_{d}$.

\subsection{Linear Foliations of Torii and Classical
Diophantine Approximation}\label{linearfoliationsoftorii}

Let $\mathcal{F}_{V}$ be the linear foliation of $\T^{p+q}$
associated to the subspace $V\subset\R^{p+q}$.  As in \S~2.1.2, we
regard $V$ as the graph of the $q\times p$ matrix ${\bf R}$.  Let
$\varphi_{\textbf{R}} : \Z^{p}\rightarrow \T^{q}$ be the homomorphism used to
define $\mathcal{F}_{V}$.  Let $\{U_{i}\}$ be a neighborhood basis
in $\T^{q}$ about $\bar{\bf 0}$.  We define a nested set
$\{ G_{i}\}\subset \Z^{p}$ by ${\bf n}\in G_{i}$ if
and only if $\varphi_{\textbf{R}} ({\bf n})\in U_{i}$.  Denote
\[ \bast\Z^{p}_{\bf R}\;\; := \;\; \bigodot G_{i}\;\;=\;\;
\bast\varphi_{\textbf{R}}^{-1}(\bast\T^{q}_{\epsilon}),  \]
a subgroup of $\bast\Z^{p}$.
If $p=q=1$ and ${\bf R}=r\in\R$, we write instead
$\bast\Z_{r}$.  

\begin{theo}\label{torusgerm} If {\rm $\textbf{R}\not\in M_{q,p}(\Q )$}, then $\fgrm
(\mathcal{F}_{V},x,\hat{x})=\bast\Z^{p}_{\bf R}$. 
Otherwise,
\[  \fgrm
(\mathcal{F}_{V},x,\hat{x})\;\; =\;\;\left\{ 
                                       \begin{array}{cl}
                                       0 & \text{if } x\not=\hat{x} \\
                                       \text{{\rm $\bast\Z^{p}_{\bf R}$}} & \text{otherwise}
                                       \end{array}
\right.  \]
\end{theo}

\begin{proof} If $\textbf{R}\not\in M_{q,p}(\Q )$, then $\varphi_{\textbf{R}}$
has dense image and the result follows by Theorem~\ref{fgrmGsusp}.
If not, then all of the leaves are torii so $\bast {\sf D}(\tilde{x},\hat{x})=0$
unless $x=\hat{x}$, in which case, if $L$ is the leaf containing $x$, 
$\bast {\sf D}(\tilde{x},\hat{x})=\bast\pi_{1}L=\bast\Z^{p}_{\bf R}$.
\end{proof}

If $\textbf{R}$ is a vector with at least one irrational entry, 
then Theorems~\ref{nonstdquot} and \ref{torusgerm} give:

\begin{coro}  Every finite dimensional torus $\T^{q}$ is algebraically isomorphic
to a quotient of the nonstardard intergers $\bast\Z$.
\end{coro}

\begin{theo}\label{ideal} {\rm ${}^{\ast}\Z^{p}_{{\textbf{R}}}$} is an ideal
in $\bast\Z^{p}$ if and only if {\rm $\textbf{R}\in M_{q,p}(\Q )$}.
\end{theo}

\begin{proof} Suppose that $\textbf{R}\in M_{q,p}(\Q )$ and let $a_{k}$ = the
l.c.d.\ of the entries of $\textbf{r}_{k}$ = the $k$th column of
$\textbf{R}$. Write
\[\mathfrak{a} = (a_{1})\oplus\cdots\oplus (a_{p})\] where
$(a_{k})$ is the ideal generated by $a_{k}$. Note that
${}^{\ast}\mathfrak{a}\subset{}^{\ast}\Z^{p}_{\textbf{R}}$.  On
the other hand, rationality of the entries of the $\textbf{r}_{k}$
implies that a sequence $\{ \textbf{n}_{\alpha}\}\subset\Z^{p}$
defines an element of ${}^{\ast}\Z^{p}_{\textbf{R}}$ if and only
if there exists $X\in\mathfrak{U}$ such that
$\varphi_{\textbf{R}}(\textbf{n}_{\alpha})=\bar{\textbf{0}}$ for all $\alpha\in
X$. This is equivalent to $\textbf{n}_{\alpha}\in \mathfrak{a}$
for all $\alpha\in X$. Thus
${}^{\ast}\Z^{p}_{\textbf{R}}={}^{\ast}\mathfrak{a}$ which is an
ideal in ${}^{\ast}\Z^{p}$.

Suppose now that $\textbf{r}=\textbf{r}_{k}\notin\Q^{q}$ for some
$k$, $1\leq k\leq p$. Let $\{ \textbf{n}_{\alpha}\}$ represent an
element $\bast\textbf{n}\in \bast\Z_{\textbf{R}}^{p}$, and denote by $\{
n_{\alpha}\}$ the sequence of $k$-th coordinates of the
$\textbf{n}_{\alpha}$ . Note that
$\overline{n_{\alpha}\textbf{r}}\not= \bar{\textbf{0}}$ for all
$\alpha$ since $\textbf{\textbf{r}}$ is not rational.  In fact,
for any $\delta>0$ we may find a sequence of integers $\{
m_{\alpha}\}$ such that
$\overline{m_{\alpha}n_{\alpha}\textbf{r}}$ is {\em not} within
$\delta$ of $\bar{\textbf{0}}$.  Let
$\textbf{m}_{\alpha}\in\Z^{p}$ be the vector whose $k$th
coordinate is $m_{\alpha}$ and whose other coordinates are $0$.
Then the sequence $\{
\textbf{m}_{\alpha}\cdot\textbf{n}_{\alpha}\}$ does not converge
with respect to $\{ G_{i}\}$ {\em i.e.} $\bast\textbf{m}\cdot\bast\textbf{n}\not\in
\bast\Z^{p}_{\textbf{R}}$, so
${}^{\ast}\Z^{p}_{\textbf{R}}$ is not an ideal.
\end{proof}

Theorem~\ref{ideal} draws another sharp distinction between
$\Z$ and ${}^{\ast}\Z$: every subgroup of the former is an ideal,
while this is false for the latter.

We spend the rest of this section studying $\bast\Z^{p}_{\bf R}$,
in and of itself a complicated and inntriguing
object.  Let us begin with the following
alternate description of ${}^{\ast}\Z^{p}_{\textbf{R}}$: 
\begin{equation}\label{altdes}
{}^{\ast}\Z^{p}_{\textbf{R}}\;\;=\;\; \Big\{{}^{\ast}\textbf{n}\in
{}^{\ast}\Z^{p}\; \Big|\;\; \exists
\;{}^{\ast}\textbf{n}^{\bot}\in{}^{\ast}\Z^{q}\text{ such that }
\textbf{R}({}^{\ast}\textbf{n})-{}^{\ast}\textbf{n}^{\bot}\in
{}^{\ast}\R^{q}_{\epsilon}\Big\}.\end{equation}  Given
${}^{\ast}\textbf{n}\in{}^{\ast}\Z^{p}_{\textbf{R}}$, the
corresponding element
${}^{\ast}\textbf{n}^{\bot}\in{}^{\ast}\Z^{q}$ is called the {\em
dual} of ${}^{\ast}\textbf{n}$; it is uniquely determined. From
(\ref{altdes}), it is clear that the set
\[ ({}^{\ast}\Z^{p}_{\textbf{R}})^{\bot}\;\; :=\;\;\Big\{
{}^{\ast}\textbf{n}^{\bot}\; \Big|\;\;
 {}^{\ast}\textbf{n}^{\bot}\text{ is the dual of }
 {}^{\ast}\textbf{n}\in {}^{\ast}\Z^{p}_{\textbf{R}}\Big\} \]
 is a subgroup of ${}^{\ast}\Z^{q}$, called the {\em dual} of
 ${}^{\ast}\Z^{p}_{\textbf{R}}$. Note that when $\textbf{R}\in M_{q,p}(\R\setminus\Q )$ has a
left-inverse $\textbf{S}$, we have
 $({}^{\ast}\Z^{p}_{\textbf{R}})^{\bot}={}^{\ast}\Z^{q}_{\textbf{S}}$.

 Similarly, the set
 \[ {}^{\ast}\R^{q}_{\textbf{R},\epsilon}\;\;=\;\;\Big\{
 {}^{\ast}\textbf{$\epsilon$}\in{}^{\ast}\R^{q}_{\epsilon} \;
 \Big|\;\;
 \exists{}^{\ast} \textbf{n}\in{}^{\ast}\Z^{p}_{\textbf{R}}\text{
 such that } \textbf{R}({}^{\ast} \textbf{n})-{}^{\ast}
 \textbf{n}^{\bot}={}^{\ast}\textbf{$\epsilon$}\Big\}\]
 is a subgroup of ${}^{\ast}\R_{\epsilon}^{q}$, called the {\em group of
rates} of $\textbf{R}$.

 The following proposition is an immediate consequence of (\ref{altdes}).

\begin{prop}  The maps {\rm ${}^{\ast} \textbf{n}\mapsto{}^{\ast}
\textbf{n}^{\bot}$} and {\rm ${}^{\ast}
\textbf{n}\mapsto{}^{\ast}\mathbf{\epsilon}$} define isomorphisms
{\rm
\[   {}^{\ast}\Z^{p}_{\textbf{R}}\;\;\cong\;\;
({}^{\ast}\Z^{p}_{\textbf{R}})^{\bot}\;\;\;\;\;\;\text{    and
}\;\;\;\;\;\;
{}^{\ast}\Z^{p}_{\textbf{R}}\;\;\cong\;\;{}^{\ast}\R^{q}_{\textbf{R},\epsilon}.\]}
\end{prop}

\begin{note}[A.Verjovsky]\label{dioph} Using formulation (\ref{altdes}) of
${}^{\ast}\Z^{p}_{\textbf{R}}$, it follows that every triple
\[({}^{\ast} \textbf{n},\, {}^{\ast} \textbf{n}^{\bot},\,
{}^{\ast}\mathbf{\epsilon})\] represents a
diophantine approximation of $ \textbf{R}$.   Thus we may regard
${}^{\ast}\Z^{p}_{\textbf{R}}$ as the {\em group of diophantine
approximations} of $ \textbf{R}$.

For example, when $p=q=1$ and $r\in\R\setminus\Q$, ${}^{\ast}n$
and ${}^{\ast}n^{\bot}$ are equivalence classes of sequences $\{
x_{\alpha}\}$ and $\{ y_{\alpha}\}$ $\subset\Z$, and
${}^{\ast}\epsilon$ an equivalence class of sequence $\{
\epsilon_{\alpha}\}\subset\R$, $\epsilon_{\alpha}\rightarrow 0$, such
that
\[ \left| r-\frac{y_{\alpha}}{x_{\alpha}}\right|\;\;=
\;\;\left|\frac{\epsilon_{\alpha}}{x_{\alpha}}\right|\;\;\longrightarrow\;\;
0.\] Conversely, every diophantine
approximation of $r$ defines uniquely a triple
$({}^{\ast}n,{}^{\ast}n^{\bot},{}^{\ast}\epsilon )$.
\end{note}

Recall that two irrational numbers $r,s\in\R\setminus\Q$ are {\em
equivalent} if there exists \[A\;\; =\;\;\left(
\begin{array}{cc}
           a & b \\
           c & d
          \end{array}
                       \right) \;\in\; {\sf SL}(2, \Z )\] such that $s=A(r)=(ar
                       +b)/(cr+d)$.

\begin{prop}\label{equiv}  If $r$ and $s$ are equivalent irrational
numbers, then
${}^{\ast}\Z_{r}\cong {}^{\ast}\Z_{s}$.
\end{prop}

\begin{proof}  Given ${}^{\ast}n\in {}^{\ast}\Z_{r}$, observe that
$(cr+d){}^{\ast}n\simeq
c{}^{\ast}n^{\bot}+d{}^{\ast}n\in {}^{\ast}\Z.$ Write
${}^{\ast}m=c{}^{\ast}n^{\bot}+d{}^{\ast}n$.  Then $\bast m\in
{}^{\ast}\Z_{s}$, since
\[s{}^{\ast}m\;\;\simeq\;\;(ar+b){}^{\ast}n\;\;\simeq\;\;a{}^{\ast}n^{\bot}+b{}^{\ast}n\;\in\;\bast\Z.\]
The
association
${}^{\ast}n\mapsto {}^{\ast}m$ defines an injective homomorphism
$\psi :\bast\Z_{r}\rightarrow\bast\Z_{s}$, with inverse defined
$\psi^{-1}(\bast m)\simeq (-cs+a)\bast m$.
\end{proof}

\begin{note}  Two irrational numbers $r,s$ are called {\em
virtually equivalent} if there exists $A\in {\sf SL}(2,\Q )$ such that
$A(r)=s$.  In this case, there exists a pair of monomorphisms
\[ \psi_{1} : \bast\Z_{r}\hookrightarrow\bast\Z_{s}\;\;\text{ and
}\;\; \psi_{2} : \bast\Z_{s}\hookrightarrow\bast\Z_{r},\] defined
as in Proposition~\ref{equiv}. In other words, $\bast\Z_{r}$ and
$\bast\Z_{s}$ are {\em virtually isomorphic}. These maps are
mutually inverse to each other if and only if $A\in {\sf SL}(2,\Z )$.
\end{note}

We are led to make the following conjecture.

\begin{conj}\label{Diophantineconjecture}  If ${}^{\ast}\Z_{r}\cong {}^{\ast}\Z_{s}$ for irrational
numbers $r$, $s$, then $r$
and $s$ are equivalent.
\end{conj}

A verified Conjecture~\ref{Diophantineconjecture} would augur
a group theoretic approach to diophantine approximation.

\subsection{Anosov Foliations and Hyperbolic Diophantine Approximation}

Let $\Gamma$ be a discrete subgroup of ${\sf PSL}(2,\R )$ with no elliptics, 
$\Sigma=\Gamma\backslash\HP^{2}$ the
corresponding Riemann surface.  Let $\rho :\Gamma\rightarrow {\sf
Homeo}(\SC^{1})$ be the representation of $\Gamma$ on
$\SC^{1}\approx \partial\HP^{2}$ and denote as in
\S~\ref{Anosov} the associated Anosov foliation by
$\mathcal{F}_{\Gamma}$. Fix $t,\xi\in \SC^{1}$, consider a
neighborhood basis $\{ U_{i}(\xi )\}$ about $\xi$, and define the nested
set $\{ G_{i}(t;\xi )\}\subset\Gamma$ by
\[  G_{i}(t; \xi )\;\; =\;\; 
\big\{ A\in\Gamma\; \big| \;\; \rho_{A}(t)\in U_{i}(\xi )\big\}.\]

\begin{prop}  Let $\hat{x}\in \mathcal{F}_{\Gamma}$ be contained in a leaf 
covered by $\HP\times \{ \xi\}$ and let $x$ be contained in a leaf covered by
by $\HP\times \{ t\}$.  Then 
\[\fgrm
(\mathcal{F},x,\hat{x})\;\;=\;\;\bigodot \big( G_{i}(t; \xi )\cdot  G_{i}(t; \xi )^{-1}  \big).\]
\end{prop}

\begin{proof}  Immediate from the definition of $\fgrm$.
\end{proof}

Classically \cite{Patterson}, given $\xi\in \SC^{1}$ in the limit set of
$\Gamma$ and $t\in \SC^{1}$, a $\Gamma$-{\em hyperbolic diophantine approximation} of
$\xi$ based at $t$ is a sequence $\{ A_{\alpha}\}\subset \Gamma$ such that
$|\xi -A_{\alpha}(t)|\rightarrow 0$, where $|\cdot |$ is the norm
induced by the inclusion $\SC^{1}\subset\R^{2}$. 
It follows from our definitions that
$\bast {\sf D}(\tilde{x}, \hat{x})$ consists precisely of equivalence classes
of $\Gamma$-hyperbolic diophantine
approximations.

\section{The Fundamental Germ of a Quasisuspension}

Let $\mathcal{L}_{\rho}$ be a Galois suspension, $\mathcal{X}\subset \mathcal{L}_{\rho}$
a $\pi_{1}B$ invariant closed set, $\mathcal{L}_{0}=\mathcal{L}_{\rho}\setminus 
\mathcal{X}$.
Let $H<{\sf Homeo}_{\omega{\sf -fib}}({\cal L}_{0})$
be a subgroup acting properly discontinuously and
let $\mathcal{Q}=H\backslash\mathcal{L}_{0}$ 
be the resulting quasisuspension.  See \S 2.2.
We have the following analogue of 
Proposition~\ref{nodependenceonT}:

\begin{prop}\label{nodependenceonT-qs}   If $T$ and $T'$ are 
diophantine transversals containing $x$ and $\hat{x}$ then
\[\fgrm (\mathcal{Q},x,\hat{x},  T)\;\; =\;\;
\fgrm (\mathcal{Q},x,\hat{x}, T').\]
\end{prop}

\begin{proof}  First suppose
that the leaf $L$ containing $x$ has the same topology as any leaf 
$L_{0}$ lying above
it in $\mathcal{L}_{0}$: in other words, $H_{L}=1$. 
Then the diophantine group $G_{L}$ is generated only by elements of 
$\pi_{1}B$ and $\pi_{1}L_{0}$. 
We may assume that the transversal $T$ lifts to 
an $H$ orbit of disjoint $\pi_{1}B$ transversals $H\cdot T_{0}$ in 
$\mathcal{L}_{0}$, 
wherein it follows
that
\begin{equation}\label{liftedtriple} 
\fgrm (\mathcal{Q},x,\hat{x},  T) \;\; = \;\;
\fgrm (\mathcal{L}_{0},x_{0},\hat{x}_{0},  T_{0})
\end{equation}
where $(x_{0},\hat{x}_{0}, T_{0})$ is a triple that covers
$(x,\hat{x}, T)$.  On the other hand, 
since the $\pi_{1}B$-invariant set $\mathcal{X}$ which we removed 
from $\mathcal{L}_{\rho}$ to get $\mathcal{L}_{0}$ is closed, we may assume
that $T_{0}$ is a diophantine transversal for $\mathcal{L}_{\rho}$.
It follows then that
\[\fgrm (\mathcal{L}_{0},x_{0},\hat{x}_{0},  T_{0})
\;\; =\;\; \bast\pi_{1}L_{0}\cdot
\fgrm (\mathcal{L}_{\rho},x,\hat{x}, T_{0}) .
 \]
The same is true for $T'$ 
so by Proposition~\ref{nodependenceonT} the result follows.
 
 Now suppose that $H_{L}\not=1$.
 Then there are $\pi_{1}B$ transversals $T_{0}, T_{0}'\subset\mathcal{L}_{0}$ 
 covering $T, T'$ such that every 
 $G_{L}$-diophantine approximation of $\hat{x}$ along $T$ resp.\ $T'$
 is of the form 
 \[\bast\gamma\cdot\bast h\cdot \bast g ,\;\; \text{resp.}\;\;
 \bast\gamma\cdot\bast h\cdot \bast g',\] 
 where 
 $\bast g\in \bast {\sf D}(\tilde{x},\hat{x}_{0},T_{0})$,
 $\bast g'\in \bast {\sf D}(\tilde{x},\hat{x}_{0},T_{0}')$ (here $\hat{x}_{0}\in
 T_{0}\cap T_{0}'$ covers $\hat{x}$),
 $\bast h\in\bast H_{L}$ and $\bast\gamma\in\bast\pi_{1}L_{0}$.  
 By the previous
 paragraph, we have 
 $\bast {\sf D}(\tilde{x},\hat{x}_{0},T_{0})=
 \bast {\sf D}(\tilde{x},\hat{x}_{0},T'_{0})$
 and the result follows.
 \end{proof}

Accordingly, we drop mention of $T$ and write 
$\fgrm (\mathcal{Q},x,\hat{x})$.

\begin{note} The proof of Proposition~\ref{nodependenceonT-qs} shows
that $\bast\pi_{1}(L)$ is a subgroup of 
$\fgrm (\mathcal{Q},x,\hat{x})$. 
In addition, there is a monomorphism
$\fgrm (\mathcal{L}_{\rho},x,\hat{x})
\hookrightarrow
\fgrm (\mathcal{Q},x,\hat{x})$, an isomorphism if $H_{L}=\{ 1\}$.
\end{note}

\subsection{Sullivan Solenoids and the 
Baumslag-Solitar Groups}

Consider the Baumslag-Solitar group
\[G_{\sf BS}\;\;=\;\; G_{\sf BS}(d)\;\; =\;\; \big\langle f,x\;:\;\; fxf^{-1}=x^{d}\big\rangle .\]
Define a nested set about 1 by
\begin{equation}\label{BSNestedSet}
  G_{i}\;\;=\;\;\Big\{ f^{m}x^{rd^{i}}\;\Big|\;\; m,r\in\Z \Big\} ,
\end{equation}
and denote 
\[ \gzg{G_{\sf BS}}\;\; :=\;\; \bigodot (G_{i}\cdot G_{i}^{-1}) . \]

\begin{theo}\label{Baum}  $\gzg{G_{\sf BS}}$ is a group.
\end{theo}

\begin{proof}   Observe by induction that in $G_{\sf BS}$,
\begin{equation}\label{bsrel}
x^{-d^{\alpha}}f\;\;=\;\; fx^{-d^{\alpha-1}} \end{equation} for
all $\alpha>0$.  To see that $\gzg{G_{\sf BS}}$ is a group, it suffices to check that
$G_{i}\cdot G_{i}^{-1}$ is a group for all $i$. Write a generic
element $g\in G_{i}\cdot G_{i}^{-1}$ in the form
$ g = f^{l}x^{rd^{i}}f^{m} $
for $l,m,r\in \Z$.  Then an element $gh^{-1}$, $g, h\in G_{i}\cdot
G_{i}^{-1}$ may be written (using (\ref{bsrel}))
\[ gh^{-1}\;\;=\;\;f^{l}x^{rd^{i}}f^{m}x^{sd^{i}}f^{n}\;\;=\;\;
\left\{
\begin{array}{ll}
f^{l}x^{(r+sd^{m})d^{i}}f^{m+n} & \text{if $m>0$} \\
& \\
f^{l+m}x^{(rd^{m}+s)d^{i}}f^{n}& \text{if $m\leq 0$}
\end{array}
\right.  ,
\]
where $l,m,n,r,s\in \Z$. It follows that $gh^{-1}\in G_{i}\cdot
G_{i}^{-1}$.
\end{proof}

\begin{note}\label{notagroupoid}  The ultrascope $\bigodot G_{i}$
is not even a groupoid as elements do not have inverses.  
Indeed, consider the sequence
$ \{ g_{\alpha}\} = \big\{ f^{-m_{\alpha}}x^{d^{\alpha}}\big\}$, 
where $m_{\alpha}>\alpha>0$,
$\alpha=1,2,\dots$. Note that $\{ g_{\alpha}\}$ defines an element
of $\bigodot G_{i}$. Using (\ref{bsrel}), we may write
the inverse sequence
\[ \{ g_{\alpha}^{-1}\}\;\;=\;\;\big\{ x^{-d^{\alpha}}f^{m_{\alpha}}\big\}\;\;=\;\;
\big\{
f^{\alpha}x^{-1}f^{m_{\alpha}-\alpha}\big\}.\] Since
$m_{\alpha}>\alpha$, we cannot use the defining relation of
$G_{\sf BS}$ to move the remaining $f^{m_{\alpha}-\alpha}$ to the
left of the $x$-term. It follows that $\{ g_{\alpha}^{-1}\}$ does
not define an element of $\bigodot G_{i}$, so the
latter does not have the structure of a groupoid.
\end{note}

\begin{theo} For all $x,\hat{x}\in \D_{f}$ with $x\in L$, 
\[  \fgrm (\widehat{\D}_{f},x,\hat{x})
\;\;\cong\;\;\left\{ 
             \begin{array}{ll}
             \gzg{G_{\sf BS}} & \text{if $L$ is an annulus} \\
             & \\
             \bast\Z_{\hat{d}} & \text{if $L$ is a disk}
             \end{array}
                        \right.\]
In either event, 
$\fgrm (\widehat{\D}_{f},x,\hat{x})$ is a group.
\end{theo}

\begin{proof}  First suppose $L$ is an annulus.  
The action of $\pi_{1}\D^{\ast}\cong\Z$ on $\widehat{\D}$
is generated by $(z,\, \hat{n})\mapsto (z,\hat{n}+1)$, where
$(z,\hat{n})\in\HP^{2}\times\widehat{\Z}_{d}$.  Then if $\gamma$
is the generator of $\pi_{1}\D^{\ast}$, we have
$\hat{f}\,\gamma\,\hat{f}^{-1}=\gamma^{\; d}$.
It follows that the diophantine group 
is isomorphic to $G_{\sf BS}$. The set of diophantine approximations
$\bast {\sf D}(\tilde{x},\hat{x})$ is equal to $\bigodot G_{i}$, where $G_{i}$
is the nested set (\ref{BSNestedSet}).  The result now follows by
definition of $\fgrm$.  If $L$ is a disk, then 
$\fgrm (\widehat{\D}_{f},x,\hat{x})=
\fgrm (\widehat{\D},x_{0},\hat{x}_{0})$ where $(x_{0},\hat{x}_{0})$
covers $(x, \hat{x})$.  By the results of \S~\ref{invlimit} 
we have $\fgrm (\widehat{\D},x_{0},\hat{x}_{0})
=\bast\Z_{\hat{\epsilon}}(d)$.
\end{proof}

The example of $\widehat{\D}_{f}$ illustrates the advantage of the ``nonabelian Grothendieck
group'' type construction used in Definition~\ref{FGDefn}: by Note~\ref{notagroupoid}, 
the naive choice ``$\fgrm =\bast {\sf D}$'' would not even have produced
a groupoid.

\subsection{Reeb Foliations}

Let $\mathcal{F}_{\sf Reeb}$ be a Reeb foliation.  The diophantine group here
is $\Z$. Recall that $L_{0}$ is the torus leaf.

\begin{theo}   For any pair $x,\hat{x}\in \mathcal{F}_{\sf Reeb}$ contained
in a diophantine transversal with $x\in L$, 
\[  \fgrm (\mathcal{F}_{\sf Reeb},x,\hat{x})
\;\;\cong\;\;\left\{ 
             \begin{array}{ll}
             \bast\Z^{2} & \text{if $x=\hat{x}\in L_{0}=L$} \\
             & \\
             \bast \Z & \text{if $\hat{x}\in L_{0}\not= L$} \\
             & \\
             0 & \text{otherwise}
             \end{array}
                        \right.\]
In every case, $ \fgrm (\mathcal{F}_{\sf Reeb},x,\hat{x})$ is a group.
\end{theo}

\begin{proof}  Suppose first that $x=\hat{x}\in L_{0}$.  Then
$\fgrm (\mathcal{F}_{\sf Reeb},x,\hat{x})=\bast\pi_{1}L_{0}=\bast\Z^{2}$.
If $x\not=\hat{x}$ and $L=L_{0}$, there is no diophantine transversal
containing the two points hence the fundamental germ is undefined.  
Now if $x\in L$, $\hat{x}\in L_{0}\not= L$ are contained in a diophantine
transversal, then 
a sequence $\{ n_{\alpha}\}$ is a diophantine approximation
if and only if it is infinite. Thus $\bast {\sf D}(\tilde{x},\hat{x})=
\bast\Z_{\infty}:= {}^{\ast}\Z\setminus {}^{\ast}\Z_{\rm
fin}$, the infinite nonstandard integers.  Then
\[\fgrm (\mathcal{F}_{\sf Reeb},x,\hat{x}) \;\;=\;\;
\bast\Z_{\infty}-\bast\Z_{\infty}
\;\;=\;\; \bast \Z.\]
If $\hat{x}\in L'\not= L_{0}$, there are no accumulations of $L$
on $L'$ so the fundamental germ is 0.
\end{proof}

Intuitively, when $\hat{x}\in L_{0}\not= L$,
$\fgrm (\mathcal{F}_{\sf Reeb},x,\hat{x})$
records the approximation by the dense leaf of the circumferential
cycle $c\subset L_{0}$ through $\hat{x}$. On the other hand, 
$\fgrm (\mathcal{F}_{\sf Reeb},x,\hat{x})$ does not predict the
meridian cycle $c'\subset L_{0}$.  Instead, $c'$ is approximated by
a sequence of inessential loops in $L$ that move off to infinity,
and such sequences are not the stuff of $\fgrm$.

\section{The Fundamental Germ of a Double Coset Foliation}

Let $\mathfrak{G}$ be a Lie group, $\mathfrak{H}<\mathfrak{G}$ a
closed subgroup, $\Gamma<\mathfrak{G}$ a discrete subgroup and
$\mathcal{F}_{\mathfrak{H},\Gamma}$ the associated double coset
foliation.
The situation is considerably more subtle due to the fact that the diophantine
group is no longer discrete.  Thus
two choices of diophantine transversal
$T_{1}$, $T_{2}$ through $x,\hat{x}$ yield distinct sets of
diophantine approximations, 
in contrast with the
case of a (quasi)suspension.  Note on the other hand that
{\em every} transversal is diophantine, since the universal covers
of the leaves
are homogeneous with respect to the left action of the diophantine
group $\widetilde{\mathfrak{H}}$. 
In fact, if $x_{1}$ and $x_{2}$ are contained
in the same leaf, then $\tilde{a}\cdot\tilde{x}_{1}=\tilde{x}_{2}$
for some $\tilde{a}\in \widetilde{\mathfrak{H}}$.  This yields a bijection
of diophantine sets
\[ \bast {\sf D}(\tilde{x}_{1},\hat{x},T_{1})\;
\longrightarrow\;\bast {\sf D}(\tilde{x}_{2},\hat{x},T_{2})\]
defined $\bast g_{1}\mapsto\bast g_{2}$ if
$\bstar g_{1}=\bstar g_{2}\cdot\tilde{a}$ in $\bstar\widetilde{\mathfrak{H}}$.
That is, the bijection is given by the equality 
$\bstar {\sf D}(\tilde{x}_{1},\hat{x},T_{1})=
\bstar {\sf D}(\tilde{x}_{2},\hat{x},T_{2})\cdot\tilde{a}$.
However, it is not clear that the following prescription for a map of
fundamental germs:
\begin{equation}\label{prescription} \bast u_{1}\mapsto \bast u_{2}\quad
\text{ iff }\quad\bast u_{1}=\bast g_{1}\bast h_{1}^{-1},\; 
\bast u_{2}=\bast g_{2}\bast h_{2}^{-1} \text{ and }
\bstar g_{1}=\bstar g_{2}\cdot\tilde{a},\;\bstar h_{1}=\bstar h_{2}\cdot\tilde{a}
\end{equation} is
well-defined since there might be, say, another representation
$\bast u_{1}=\bast g_{1}'(\bast h_{1}')^{-1}$ which leads to a different assignment.
Even if (\ref{prescription}) were well-defined, 
there is no reason to expect that it should respect
the groupoid structure.  When $\bstar\widetilde{\mathfrak{H}}$ is a group,
one can say more:

\begin{lemm}  If $\bstar\widetilde{\mathfrak{H}}$ is a group then
$\bast u\circ\bast v=\bast w$ in 
$\fgrm (\mathcal{F}_{\mathfrak{H},\Gamma},x,\hat{x}, T)$
implies
$\bstar u\cdot\bstar v=\bstar w$ in $\bstar\widetilde{\mathfrak{H}}$.
\end{lemm}

\begin{proof} This follows immediately since the groupoid structure
of the fundamental germ is defined in terms of left multiplication 
on the unit space $\bstar {\sf D}(\tilde{x},\hat{x}, T)$.  
\end{proof}

\begin{prop} If $\bstar\widetilde{\mathfrak{H}}$ is a group
and $T_{1}$ and $T_{2}$ are diophantine transversals through $x_{1},\hat{x}$ 
and $x_{2},\hat{x}$, respectively, where $x_{1},x_{2}$ belong to the
same leaf $L$, then
\[\fgrm (\mathcal{F}_{\mathfrak{H},\Gamma},x_{1},\hat{x}, T_{1})\;\; \cong\;\;
\fgrm (\mathcal{F}_{\mathfrak{H},\Gamma},x_{2},\hat{x}, T_{2}).\]
\end{prop}

\begin{proof}  It is clear now that the bijection (\ref{prescription}) 
is well-defined: in fact,
since $\bstar\widetilde{\mathfrak{H}}$ is a group, we have
$\bstar u_{1}=\bstar u_{2}$.  From this it follows that 
${\sf Dom} (\bast u_{1} )={\sf Dom} (\bast u_{2} )\cdot\tilde{a}$, 
and that the bijection
(\ref{prescription}) defines a groupoid isomorphism.
\end{proof}

We shall assume from this moment on that 
$\bstar\widetilde{\mathfrak{H}}$ is a group.  
We will then not mention the base point $x$ and the
transversal $T$ and write 
$\fgrm (\mathcal{F}_{\mathfrak{H},\Gamma},L,\hat{x})$
where $L$ is the leaf along which diophantine approximations are taking place.
If $\hat{x}\in L$ we write simply 
$\fgrm (\mathcal{F}_{\mathfrak{H},\Gamma},L)$.

We now give a ``diophantine'' description of 
$\bast {\sf D}(\tilde{x},\hat{x}, T)$,
similar in spirit to that of $\bast\Z^{p}_{\bf R}$ appearing in (\ref{altdes}).  
Denote by $p:\widetilde{\mathfrak{H}}\rightarrow
\mathfrak{H}$ the universal cover of $\mathfrak{H}$.  Suppose
that $L$ is covered by a coset $\mathfrak{H}g$ and
$\hat{g}\in\mathfrak{G}$ is an element covering $\hat{x}$.  A subset 
$\mathcal{T}^{\hat{g}}\subset\mathfrak{G}$ is called a {\em local section} at $\hat{g}$
for the quotient map
$\mathfrak{G}\rightarrow\mathfrak{H}\backslash\mathfrak{G}$ if 
$\mathcal{T}^{\hat{g}}$ maps
homeomorphically onto an open subset containing $\mathfrak{H}\hat{g}$.  We may
assume without loss of generality that the transversal $T$ 
through $\hat{x}$ lifts to a local section $\mathcal{T}^{\hat{g}}$
through $\hat{g}$. 
As our interest is in sequences which converge to 
$\hat{g}$ in $\mathcal{T}^{\hat{g}}$,
we may assume also that $\mathcal{T}^{\hat{g}}=\hat{g}\cdot\mathcal{T}$ for
some local section $\mathcal{T}$ about $1$.  Let 
$\bast\mathcal{T}_{\epsilon}\subset\bast\mathfrak{G}_{\epsilon}$
denote the set of infinitesimals which are represented by sequences
in $\mathcal{T}$.

Now let $\bast\tilde{h}$ be a diophantine approximation 
of $\hat{x}$ based at $\tilde{x}$ along $T$, which is
characterized by the property that $\{ p(\bast\tilde{h})\cdot g \}$ lies in
$\hat{g}\cdot \bast\mathcal{T}_{\epsilon}\cdot\bast\Gamma$.  This
gives the following diophantine description of $\bast{\sf D}(\tilde{x},\hat{x},T)$: 
\begin{equation}\label{diophexpdoubcos}\bast{\sf D}(\tilde{x},\hat{x},T)
\;\; =\;\; \left\{ \bast\tilde{h}\in \bast\widetilde{\mathfrak{H}} \;
\left|\;\;\exists\,\bast\gamma\in\bast\Gamma ,\; \bast \epsilon\in 
\bast\mathcal{T}_{\epsilon}
 \text{ such that } 
 \hat{g}^{-1}\cdot p(\bast\tilde{h}) \cdot g\cdot \bast \gamma  =
\bast \epsilon
\right.\right\}.
\end{equation}

The element
$\bast\tilde{h}^{\perp}:=\bast\gamma$ associated
to $\bast\tilde{h}$ in (\ref{diophexpdoubcos}) is called the {\em dual}
of $\bast\tilde{h}$.
When $\hat{g}=g$, we let
$\bast\widetilde{\mathfrak{H}}_{g}:=\bast{\sf D}(\tilde{x},T)$ denote
the set of diophantine approximations and let
$\bast\widetilde{\mathfrak{H}}^{\perp}_{g}$ denote the
set of duals.  Thus if
$\sigma_{g}$ denotes the conjugation map $a\mapsto g^{-1}ag$, 
\begin{equation}\label{diophexpdoubcos1}
\bast\widetilde{\mathfrak{H}}_{g}
\;\; =\;\; \left\{ \bast\tilde{h}\in \bast\widetilde{\mathfrak{H}} \;
\left|\;\;\exists\,\bast\gamma\in\bast\Gamma ,\; \bast \epsilon\in 
\bast\mathcal{T}_{\epsilon}
 \text{ such that } 
  \sigma_{g}( p(\bast\tilde{h})) \cdot \bast \gamma =
\bast \epsilon
\right.\right\} .
\end{equation} 
In general, whether $g=\hat{g}$ or not, it follows that
\[\fgrm (\mathcal{F}_{\mathfrak{H},\Gamma},L,\hat{x})
\;\; \subset \;\;
\left\{ \bast\tilde{u}\in \bast\widetilde{\mathfrak{H}}\; \left|\;\;
\exists\bast\gamma ,\bast\eta\in\bast\Gamma ,\; 
\bast \omega\in \gzg{\mathcal{T}}
\text{ s.t.\ } 
\bast\gamma\cdot\sigma_{g}(p(\bast\tilde{u}))\cdot\bast\eta =
\bast\omega\right.\right\} , \]
where 
$\gzg{\mathcal{T}}=\bast\mathcal{T}_{\epsilon}\cdot\bast\mathcal{T}_{\epsilon}^{-1}$.
The inclusion is in general strict as the following example
shows:

\begin{exam} Consider the double
coset foliation ${\sf Geod}_{\Gamma}$, which possesses a noncompact leaf $L$
and
a pair of cycles $c_{-}$, $c_{+}$ such that $L$ coils about $c_{-}$ (about $c_{+}$)
as one goes to negative (positive) infinity in $L$,
and has no other accumulations. If 
$\hat{x}$ belongs to either $c_{-}$ or $c_{+}$, 
we have
\[ \fgrm ({\sf Geod}_{\Gamma},L,\hat{x})\;\;\cong\;\;\bast \Z,\] but if $\hat{x}\in L$,
we have
$\fgrm ({\sf Geod}_{\Gamma},L,\hat{x})=0$.  
\end{exam}

\begin{note} Since $p^{-1}(e)\cong\pi_{1}\mathfrak{H}$, we
have $\bast\pi_{1}\mathfrak{H}<  
\fgrm (\mathcal{F}_{\mathfrak{H},\Gamma}, L,\hat{x})$.
\end{note}

One can
understand the description of $\bast\widetilde{\mathfrak{H}}_{g}$
appearing in (\ref{diophexpdoubcos1}) as a nonlinear version of
(\ref{altdes}). In fact, if $\mathfrak{G}$ is a linear group of $p\times p$
matrices and $g\in\mathfrak{G}$, 
then one can think of $\bast \Z^{p}_{g}$ as defined in (\ref{altdes}) as the set
of linear diophantine approximations of $g$ (approximations of $g$ by pairs
of vectors with
respect to linear algebra),
whereas $\bast\widetilde{\mathfrak{H}}_{g}$ can be thought of as
a set of nonlinear diophantine approximations of $g$ (approximations 
of $g$ by pairs of matrices with respect to matrix algebra).

We now consider the horocyclic and geodesic flows
on the unit tanget bundle
of a riemannian surface, which are,  
as is widely appreciated, deep mathematical objects.  
It should come as no suprise that this deepness
is reflected in
their fundamental germs, which
present the most complex and intractable diophantine algebra we
have encountered thus far.
In the remainder of this section, we will attempt 
to give the reader a feel for the complexity of these fundamental germs
by walking through a sample calculation.

We restrict to the case
$\mathfrak{G}=SL(2,\R )$ and $\Gamma= SL(2,\Z )$.
See \S~\ref{DoubleCoset} for the relevant notation. Consider first
the case of the (positive) horocyclic flow ${\sf Hor}={\sf Hor}^{+}_{{\sf SL}(2,\Z)}$,
that is, $\mathfrak{H}=H=H^{+}$.
If $D$ is the subgroup
of matrices of the form
\[  \left(
\begin{array}{cc}
e^{s/2} & 0 \\ 
t & e^{-s/2}
\end{array}\right)\]
$s,t\in\R$,
then $D$ defines a local section about $1$ so we take $\mathcal{T}=D$.  
Finally, since $H\cong (\R ,+)$, we shall
simplify notation by identifying $r$ with the matrix $A_{r}$ and
write $\bast\R_{g}=\bast H_{g}$ for the set of diophantine
approximations.

Let us consider the relatively simple choice
\[ g=\left(
\begin{array}{cc}
\sqrt{2} & 1 \\
1 & \sqrt{2}
\end{array}
\right).\] The right coset of $g$ is
\[ Hg\;\;=\;\;\left\{ \left(
\begin{array}{cc}
 r+\sqrt{2}& \sqrt{2}r+1 \\
1 & \sqrt{2}
\end{array}
\right) \;\Bigg|\;\; r\in \R  \right\} .\] Since $Hg$ does not define a cycle in 
$SL(2,\R )/SL(2,\Z )$ it must be dense by a theorem of
Hedlund \cite{Hedlund}, so we can expect from $g$ 
a nontrivial set of diophantine approximations.  The
conjugate of $H$ by $g$ is
\[\sigma_{g}(H)\;\;=\;\; \left\{\left(
\begin{array}{cc}
1+\sqrt{2}r& 2r \\
-r & 1-\sqrt{2}r 
\end{array}
\right) \;\Bigg|\;\; r\in \R  \right\}.\]

In order to characterize the elements of $\bast\R_{g}$,
we shall need the following generalization of $\bast\Z_{r}$.  
Let $\OI$ be the ring of integers of a number field.  For $\bast r\in\bast\R$, 
define
\[  \bast \OI_{\bast r}\;\; =\;\;\{ \bast n\in\bast \OI\; |\;\;
\exists\;\bast n^{\perp}\in\bast \OI\text{ such that } 
\bast r\cdot\bast n-\bast n^{\perp}\in\bast\R_{\epsilon}  \}  . \]
Clearly $\bast \OI_{\bast r}$ is a subgroup of $\bast\OI$.

\begin{theo}\label{fgexpdoublecoset} Let $\OI$ be the ring of integers in $\Q (\sqrt{2})$.
Then $\bast r\in\bast \R_{g}$ if and only if there exists
$\bast\gamma=\left(
\begin{array}{cc}
{}^{\ast}a & ^{\ast}b \\
{}^{\ast}c & {}^{\ast}d
\end{array}
\right)\in {\sf SL}(2,\bast\Z)$ for which
\begin{itemize}
\item $\sqrt{2}\bast a + 2\bast c,\; \sqrt{2}\bast b + 2\bast d\in  
\bast\OI_{\bast r}\;$ 
and $\; (\sqrt{2}\bast a + 2\bast c)^{\perp} =  1-\bast a$,\\
$(\sqrt{2}\bast b + 2\bast d)^{\perp}=-\bast b$.
\[
\quad\] 
\item  $\bast c, \bast d\in\bast\Z_{\sqrt{2}}$ and $\bast c^{\perp}=1-\bast a$, 
$\bast d^{\perp}=1-\bast b$.
\[
\quad\] 
\item $\bast b=-(\sqrt{2}\bast b+2\bast d )\bast r$ and 
$(\bast a+(\sqrt{2}\bast a + 2\bast c)\bast r)
(\bast d-(\bast b+\sqrt{2}\bast d)\bast r)=1$.
\end{itemize}
\end{theo}

\begin{proof}  From (\ref{diophexpdoubcos1}), $\bast r\in
\R_{g}$ if and only if there exists $\bast\gamma\in \bast \Gamma$ and
 $\bast\epsilon ,\bast \delta\in\bast\R_{\epsilon}$ with
 \[  \left(
\begin{array}{cc}
\bast a (1+\sqrt{2}\bast r) +2\bast c \bast r & \bast b (1+\sqrt{2}\bast r) +2\bast d \bast r \\
& \\
-\bast a \bast r + \bast c (1-\sqrt{2}\bast r)  & -\bast b \bast r + \bast d (1-\sqrt{2}\bast r)
\end{array}
\right)  \;\; =\;\; 
\left(
\begin{array}{cc}
1+\bast \epsilon & 0 \\
& \\
\bast\delta  & (1+\bast \epsilon)^{-1}
\end{array}
\right) .
\]
The first and third items follow immediately.  The second item follows
upon noting that we may eliminate $\bast r$ by multiplying the second
row equations by $\sqrt{2}$ and adding them
to the first row equations.
\end{proof}

Theorem~\ref{fgexpdoublecoset} illustrates why it is so difficult to say
anything about the algebraic structure of
$\bast \R_{g}$ or $\fgrm ({\sf Hor},L)$.   In order
to determine whether the sum $\bast r +\bast s$
defines an element of $\bast \R_{g}$, 
we must find a way to ``compose'' 
the corresponding duals $\bast r^{\perp} ,\bast s^{\perp}\in \R_{g}^{\perp}$ 
to obtain one for their sum, and it is not even clear
what this operation on matrices should be.  One could reverse the logic
and ask if the product $\bast r^{\perp}\cdot\bast s^{\perp}$
defines an element of $\R_{g}^{\perp}$, however this seems just
as hopeless since the diophantine conditions spelled out in the statement of 
Theorem~\ref{fgexpdoublecoset} are not stable with respect to 
matrix multiplication.

As for the geodesic flow, we leave it to the reader to
formulate the appropriate analogue of 
Theorem~\ref{fgexpdoublecoset} {\em e.g.}
using the local section $\mathcal{T}$ for which
\[\bast\mathcal{T}_{\epsilon}\;\; =\;\; 
\left\{\left(
\begin{array}{cc}
1 & \bast\delta \\
\bast\delta' & 1+\bast\delta\bast\delta'
\end{array}
\right) \;\Bigg|\;\; \bast\delta ,\bast\delta'\in \bast\R_{\epsilon}  \right\} .
\]
The result would be a set of diophantine conditions at least as daunting as that obtained for the horocyclic
flow.

\section{The Fundamental Germ of a Locally Free Lie Group Action}\label{LieAction}

The discussion here is very similar to that for a double coset,
so we will be brief.
Let $\mathfrak{B}$ be a Lie group of dimension $k$, $M^{n}$ an
$n$-manifold, $n>k$, $X\subset M^{n}$.  Let $\theta
:\mathfrak{B}\rightarrow {\sf Homeo}(X)$ be a locally-free
representation whose orbits either coincide or are disjoint
and let $\mathcal{L}_{\mathfrak{B}}$ be the
associated lamination on $X$.  Any diophantine transversal through
$x,\hat{x}$ may be
obtained as the intersection of $\mathcal{L}_{\mathfrak{B}}$ with 
a submanifold $T$ of $M^{n}$
of dimension $n-k$ such that $x,\hat{x}\in T$ and 
$T\cap (\theta (\mathfrak{B})\cdot x)$ is
discrete in $\theta (\mathfrak{B})\cdot x$.  As in the case of a
double coset foliation, when $\bstar\mathfrak{B}$ is group,
\begin{enumerate}
\item Groupoid multiplication in the fundamental germ corresponds
to multiplication in $\bstar\mathfrak{B}$. 
\item If $T_{1}$, $T_{2}$ are transversals through
$x_{1},\hat{x}$ and $x_{2},\hat{x}$ where $x_{1},x_{2}$
belong to the same leaf $L$ then 
\[\fgrm (\mathcal{L}_{\mathfrak{B}},x_{1},\hat{x}, T_{1})\;\; \cong\;\;
\fgrm (\mathcal{L}_{\mathfrak{B}},x_{2},\hat{x}, T_{2}).\]
\end{enumerate}
Accordingly we shorten to 
$\fgrm (\mathcal{L}_{\mathfrak{B}},L,\hat{x})$.
\begin{theo}  Let $\Sigma=\Gamma\backslash\HP^{2}$ be a compact hyperbolic surface, 
$\mathfrak{l}\subset\Sigma$ a geodesic lamination, $\hat{x}\in\mathfrak{l}$
and $l\subset \mathfrak{l}$ a leaf.  Then 
\[ \fgrm (\mathfrak{l},l,\hat{x})\;\; =\;\; 
\fgrm ({\sf Geod}_{\Gamma},L, \hat{v}) \]
where 
$L$ is a leaf covering $l$ and $\hat{v}$ is a tangent vector to $l$
at $\hat{x}$.
\end{theo}

\begin{proof}  This follows immediately from the fact
that any diophantine approximation of $\hat{v}$ along $L$
canonically defines a diophantine approximation of $\hat{x}$
along $l$ and {\em vice verca}.
\end{proof}

\section{Functoriality}

We begin by recalling the notion of morphism in the category of laminations. 
A {\em lamination map} $F:\mathcal{L}\rightarrow\mathcal{L}'$
is a map satisfying the following conditions:
\begin{enumerate} 
\item For every leaf $L\subset\mathcal{L}$, there exists
a leaf $L'\subset\mathcal{L}'$ with
$F(L)\subset L'$.  
\item  For
all $x\in\mathcal{L}$, there exist open transversals $T\ni
x$, $T'\ni F(x)$, such that $F(T)\subset T'$.  
\end{enumerate}

The projection $P:\mathcal{L}\rightarrow B$ of a
suspension onto its base is a lamination map.  On the other
hand, let $\mathcal{F}$ be a foliation, $M$ the underlying
manifold. Then the
canonical inclusion $\imath :\mathcal{F}\rightarrow M$ is a
map which maps leaves into the unique leaf $M$, 
yet is not a lamination map since no open transversal of
$\mathcal{F}$ is mapped into a point, an open transversal of
$M$.

Let \[F:(\mathcal{L},x,\hat{x})\longrightarrow(\mathcal{L}',x',\hat{x}')\] be a
lamination map.  We say that $F$ is {\em diophantine} if there exist
diophantine transversals $T\ni x,\hat {x}$ and  $T'\ni x',\hat{x}'$ 
such that $F(T)\subset T'$.  Note that this condition is always
satisfied if either $\mathcal{L}$ or $\mathcal{L}'$ are laminations
defined by double cosets or locally free Lie group actions.  
Denote by $L$ and $L'$ the leaves containing $x,x'$ and let
$\widetilde{F}:\widetilde{L}\rightarrow\widetilde{L}'$ be the lift
of the restriction $F|_{L}$.  Let
$\widetilde{T}\subset\widetilde{L}$,
$\widetilde{T}'\subset\widetilde{L}'$ be the pre-images of
$T\cap L$, $T'\cap L'$.
Then for $F$ diophantine there is a well-defined map
\[\bast {\sf D}F:
\bast {\sf D}(\tilde{x},\hat{x}, T )\longrightarrow 
\bast {\sf D}(\tilde{x}',\hat{x}',T' ) \]
of diophantine approximations.  If the assigment
\[ \bast u=\bast g\cdot\bast h^{-1}\;\longmapsto \; 
\bast {\sf D}F(\bast g)\cdot (\bast {\sf D}F(\bast h))^{-1} \] leads to a well-defined map
\[\gzg{F}:\fgrm (\mathcal{L},x,\hat{x})\longrightarrow
\fgrm (\mathcal{L}',x',\hat{x}'),\]
we say that $F$ is {\em germ}.

\begin{prop}\label{sushom}  Let $\mathcal{L}=\widetilde{B}\times_{\rho} F$ be a
suspension with $x,\hat{x}$ lying over $x_{0}\in B$.  Then the projection
$\xi :(\mathcal{L},x,\hat{x})\rightarrow (B,x_{0})$ is germ, and the induced
map $\gzg{\xi}$ is a groupoid
monomorphism.
\end{prop}

\begin{proof} It is clear from the definitions that 
$\bast {\sf D}\xi$ is the inclusion
\[ \bast {\sf D}(\tilde{x},\hat{x})\;\;\subset\;\; \bast\pi_{1}(B,x). \]
In particular, it follows that $\gzg{\xi}$ is well-defined.
Since the product in $\fgrm (\mathcal{L},\hat{x}, L)$ is induced by
multiplication in $ \bast\pi_{1}(B,x)$,  $\gzg{\xi}$ is a
groupoid homomorphism as well.
\end{proof}

Unfortunately, we cannot assert in general that
the map $\gzg{F}$ induced by a germ lamination map $F$ defines a
groupoid homomorphism. We now introduce a class of
lamination maps which is sufficiently well-behaved so as to allow
us to say more.

Let $\mathcal{F}$ be a foliation, $M$ the underlying space of
$\mathcal{F}$, and $\imath :\mathcal{F}\rightarrow M$ the
inclusion.  Although $\imath$ is not a lamination map, we may
nevertheless define a map of diophantine approximations as follows.   
An element $\bast g\in 
\bast {\sf D}(\tilde{x},\hat{x},T)$, represented say by $\{ g_{\alpha}\}$,
may be regarded as made up from an equivalence class of sequence 
$\{\gamma_{g_{\alpha}}\}$ where the $\gamma_{g_{\alpha}}$
are homotopy classes of curves
lying within $L$ whose endpoints converge to $\hat{x}$.
One may assume that
there is an open disc
$O\subset M$ about $\hat{x}$ such that the endpoints of
these sequences lie entirely in $O$.
By connecting their endpoints
to $\hat{x}$ by a paths contained in $O$, we obtain a sequence of
homotopy classes of curves $\{ \eta_{g_{\alpha}}\}\subset\Pi_{1} (M,x,\hat{x})$ = 
the set of homotopy classes of paths from $x$ and $\hat{x}$, 
hence a map
\[ \bast{\sf D}\imath :\bast {\sf D}(\tilde{x},\hat{x},T)\;\longrightarrow\;
\bast\Pi_{1} (M,x,\hat{x}) ,\quad\quad\quad\quad \bast g\longmapsto
\eta_{\bast g}\] 
which depends neither on $O$ nor on the
choice of connecting paths.
More generally, given $\mathcal{L}$ a lamination and $\imath
:\mathcal{L}\rightarrow X$ a map into a path-connected space, we
may define a map
$\bast {\sf D}\imath:\bast {\sf D}(\tilde{x},\hat{x},T)\rightarrow
\bast\Pi_{1}(X,\imath (x),\imath (\hat{x}))$. We say that the map
$\imath$ is {\em germ} if $\bast {\sf D}\imath$ induces a well-defined map
\[ \gzg{\imath}:\fgrm (\mathcal{L},x,\hat{x},T)\rightarrow 
\bast\pi_{1}(X,x) ,\quad\quad\quad\quad \bast u = \bast g\bast h^{-1}\longmapsto
\bast {\sf D}\imath (\bast g)\cdot (\bast {\sf D}\imath (\bast h))^{-1}.\]

\begin{defi}  Let $\mathcal{L}$ be a lamination arising from a group action, 
$X$ a path connected space.
A map $\imath :(\mathcal{L},x,\hat{x})\rightarrow (X,\imath (x),\imath (\hat{x}))$ 
is called a {\bf
fidelity} if it is germ and $\gzg{\imath}$ is a groupoid monomorphism.  We say
that $\mathcal{L}$ is {\bf faithful} if it has a fidelity.
\end{defi}

For example, by Proposition~\ref{sushom} any suspension is faithful,
however if the underlying space of a suspension $\mathcal{L}$ is a manifold $M$, 
we shall see that it is much more useful to be able to assert that the inclusion 
$\mathcal{L}\hookrightarrow M$
is a fidelity.

For the remainder of the section, the base points $x$ and $\hat{x}$ will be supressed
in order to simplify notation.

\begin{prop}  Let $\mathcal{F}_{V}$ be the foliation of $\T^{p+q}$
induced by the $p$-plane $V\subset \R^{p+q}$.  Then the inclusion
$\imath :\mathcal{F}_{V}\rightarrow\T^{p+q}$ is a fidelity.
\end{prop}
\begin{proof}  Recall that for some $q\times p$ matrix $\bf{R}$, 
$\fgrm (\mathcal{F}_{V})=\bast\Z^{p}_{\bf R}$.
Then for $\bast {\bf n}\in\bast\Z^{p}_{\bf R}$, the map
$\gzg{\imath}$ is
\[\gzg{\imath}\big({}^{\ast}{\bf
n}\big)\;\;=\;\;\big(\bast{\bf n},{}^{\ast}{\bf
n}^{\bot}\big)\;\;\in\;\;\bast\Z^{p+q}\;\;=\;\;\bast\pi_{1}\T^{p+q}
,\] where ${}^{\ast}{\bf n}^{\bot}$ is the dual to ${}^{\ast}{\bf
n}$. $\gzg{\imath}$ is then clearly an injective homomorphism.
\end{proof}

The problem of the existence of fidelities for laminations arising
from group actions is interesting but seems difficult. 

\begin{conj}  Every lamination arising form a group action 
is faithful.
\end{conj}

\begin{defi} A germ lamination map 
$F:\mathcal{L}\rightarrow \mathcal{L}'$ is
{\bf trained} if $\mathcal{L}$ and $\mathcal{L}'$ are faithful,
and there exist fidelities $\imath:\mathcal{L}\rightarrow X$,
$\imath' :\mathcal{L}'\rightarrow X'$ and a map $f:X\rightarrow
X'$ such that
\begin{equation}\label{comp}
{}^{\ast}f\circ \gzg{\imath}\;\; =\;\;\gzg{\imath'}\circ
\gzg{F}.\end{equation} The triple $(\imath ,\imath ',f)$ is called
a {\bf training} for $F$.
\end{defi}

\begin{theo}\label{funcstab}  Let 
$F:\mathcal{L}\rightarrow \mathcal{L}'$ be a
trained lamination map.  Then the induced map $\gzg{F}$ is a
groupoid homomorphism.
\end{theo}

\begin{proof}  Let $(\imath ,\imath ',f)$ be a training for $F$.
Then for all $\bast u ,\bast v\in\fgrm (\mathcal{L})$ such that
$\bast u\cdot\bast v$ is defined we have
\[ \gzg{\imath'}\circ\gzg{F}\Big(\bast u\cdot\bast v\Big)  \;\; =\;\;
\gzg{\imath'}\Big(\gzg{F}\bast u\cdot\gzg{F}\bast v\Big).\]
Since $\gzg{\imath'}$ is injective, 
$\gzg{F}\big(\bast u\cdot\bast v\big)=\gzg{F}\bast u\cdot\gzg{F}\bast v$.
\end{proof}

\begin{coro}  Let $F:(\mathcal{F},x)\rightarrow(\mathcal{F}',x')$ be a map
of foliations.  Suppose that the inclusions into the underlying
manifolds $\imath :\mathcal{F}\rightarrow M$, $\imath'
:\mathcal{F}'\rightarrow M'$ are fidelities.  Then $\gzg{F}$ is a
groupoid homomorphism.
\end{coro}

\begin{proof}  Take $f:M\rightarrow M'$ to be $F$, viewed as a map
on underlying manifolds.  Then $(\imath ,\imath ',f)$ is a
training.
\end{proof}

\begin{coro} Any map
$F:\mathcal{F}_{V}\rightarrow\mathcal{F}_{V'}$ of linear foliations
of torii induces a homomorphism $\gzg{F}$ of fundamental germs.
\end{coro}

\section{The Germ Universal Cover}

We assume throughout this section that
\begin{enumerate}
\item $\mathcal{L}$ is a weakly-minimal lamination arising from a group
action.
\item $x=\hat{x}\in L$ a fixed dense leaf.  
\end{enumerate}
We abreviate the associated fundamental germ to $\fgrm (\mathcal{L})$.
An ultrafilter $\mathfrak{U}$ is fixed throughout.

Let $p:\widetilde{L}\rightarrow L$ be the universal cover.  A sequence 
$\{ \tilde{x}_{\alpha}\}\subset\widetilde{L}$ is
called
{\em $\mathcal{L}$-convergent} if 
it projects to a sequence in $L$ converging to some 
$\hat{x}\in\mathcal{L}$.
Two $\mathcal{L}$-convergent sequences $\{ \tilde{x}_{\alpha}\}$ and 
$\{ \tilde{x}_{\alpha}'\}\subset\widetilde{L}$ are called
$\mathcal{L}$-{\em asymptotic} if their projections converge to the same point $\hat{x}$
and if for every flowbox $O$ in $\mathcal{L}$ about $\hat{x}$, there exists
$X\in\mathfrak{U}$ such that $\tilde{x}_{\alpha}$ and $\tilde{x}_{\alpha}'$
lie in a common lift of a plaque of $O$, for all $\alpha\in X$.
The asymptotic class corresponding to $\{ \tilde{x}_{\alpha}\}$ is denoted
$\bstar\tilde{x}$; we refer to $\hat{x}$ as the {\em limit} of $\bstar\tilde{x}$
and write $\lim \bstar\tilde{x}=\hat{x}$.  The set of $\bstar\tilde{x}$
with limit $\hat{x}$ is denoted ${\sf Lim}_{\hat{x}}$.

\begin{defi}  The {\bf germ universal cover} of $\mathcal{L}$ with respect
to $L$ is
\[  \gcom{\mathcal{L}}\;\; =\;\; \Big\{ \text{classes }\bstar\tilde{x}\text{ of }
\mathcal{L}\text{-convergent sequences in }\widetilde{L} \Big\} . \]
\end{defi}

Note that for any $\hat{x}\in\mathcal{L}$, every $G_{L}$-diophantine approximation 
$\bast g$ of $\hat{x}$ determines
an element of $\gcom{\mathcal{L}}$, and the sets 
${\sf Lim}_{\hat{x}}$ and $\bast {\sf D}(\tilde{x},\hat{x},T)$ are in bijective correspondence,
for any diophantine transversal $T$ through
$x,\hat{x}$.  

\begin{prop}  Let $\mathcal{L}$ be compact and suppose
that $L=\mathfrak{G}$ is a topological group for which $\bast \tilde{b},
\bast\tilde{c}\in\bast\widetilde{\mathfrak{G}}$ are $\mathcal{L}$-asymptotic 
if and only if $\bast\tilde{b}\cdot
\bast\tilde{c}^{-1}\in\bast \widetilde{\mathfrak{G}}_{\epsilon}$. Then
$\gcom{\mathcal{L}}=\bstar\widetilde{\mathfrak{G}}$.
\end{prop}

\begin{proof}  Suppose that there is some
$\bstar\tilde{b}\in \bstar\widetilde{\mathfrak{G}}$ represented
by a sequence $\{\tilde{b}_{\alpha}\}$ which is not $\mathcal{L}$-convergent.
Thus if $\{b_{\alpha}\}$ is the projection of this sequence to 
$\mathfrak{G}\subset\mathcal{L}$,
then for all $\hat{x}\in\mathcal{L}$, $\hat{x}$ has a neighborhood 
$U_{\hat{x}}\subset\mathcal{L}$
for which there is no $X\in\mathfrak{U}$ with 
$\{ b_{\alpha}\}|_{X}\subset U_{\hat{x}}$.
The $U_{\hat{x}}$ cover $\mathcal{L}$ so that there is a subcover
$U_{\hat{x}_{1}},\dots ,U_{\hat{x}_{n}}$; this implies
that there exists a partition $X_{1}\sqcup\dots\sqcup X_{n}$ of $\N$ with
$\{ b_{\alpha}\}|_{X_{i}}\subset U_{\hat{x}_{i}}$.  Since $\mathfrak{U}$
is an ultrafilter, one of the $X_{i}$ belongs to $\mathfrak{U}$, contradiction.
Thus every element $\bstar \tilde{b}\in\bstar\widetilde{\mathfrak{G}}$ defines
an element of $\gcom{\mathcal{L}}$.  Since the relation of being $\mathcal{L}$-asymptotic
coincides with differing by an infinitesimal, we are done.
\end{proof}

For example, if $\mathcal{F}_{V}$ is a linear $n$-foliation of a torus,
$\gcom{\mathcal{F}_{V}}=\bstar\R^{n}$.

Denote by
\[ \bstar p:\gcom{\mathcal{L}}\;\longrightarrow\;\mathcal{L}\]
the natural projection defined $\bstar\tilde{x}\mapsto \lim\bstar\tilde{x}$. 
The {\em leaf} $L_{\bstar\tilde{x}}$ through $\bstar\tilde{x}$
is defined to be the set of 
$\bstar\tilde{y}$ such that
\begin{enumerate}
\item If $\hat{x}=\lim \bstar\tilde{x}$ and 
$\hat{y}=\lim \bstar\tilde{y}$ then $L_{\hat{x}}=L_{\hat{y}}$.
\item There are
representative sequences $\{ \tilde{x}_{\alpha}\}$, $\{ \tilde{y}_{\alpha}\}$,
and paths $\tilde{\eta}_{\alpha}$ connecting $\tilde{x}_{\alpha}$ to 
$\tilde{y}_{\alpha}$ so that $p(\tilde{\eta}_{\alpha})$ converges
to a path connecting $\hat{x}$ to $\hat{y}$.
\end{enumerate}

\begin{theo}\label{guctopology}  $\gcom{\mathcal{L}}$ may be given the
structure of a lamination whose leaves
are nowhere dense and for which 
$\bstar p$ is an open lamination map.  
\end{theo}

\begin{proof}  Denote by $\gzg{T}\subset \gcom{\mathcal{L}}$ the pre-image of a 
transversal $T\subset\mathcal{L}$ and well-order each ${\sf Lim}_{\hat{x}}$
for $\hat{x}\in T$.  Note that the cardinalities
of the ${\sf Lim}_{\hat{x}}$ are the same: that of the continuum, 
since $L$ is dense and $L\cap T$ is countable.  We define a decomposition
\begin{equation}\label{decomposition} \gzg{T}\;\; =\;\; \bigsqcup T_{\alpha}  \end{equation}
where $T_{\alpha}$ is the section over $T$ defined by
$\hat{x}\mapsto$ the $\alpha$th element of ${\sf Lim}_{\hat{x}}$. 
 By definition of the leaves
of $\gcom{\mathcal{L}}$, given $\hat{x},\hat{y}\in T$,
\begin{equation}\label{leaves}
\left(\bigcup_{\bstar\tilde{x}\in {\sf Lim}_{\hat{x}}}L_{\bstar\tilde{x}}\right)
 \;\bigcap\;
\left(\bigcup_{\bstar\tilde{y}\in {\sf Lim}_{\hat{y}}}L_{\bstar\tilde{y}}\right)
\;\;\not=\;\;\emptyset 
\end{equation}
if and only if $L_{\hat{x}}=L_{\hat{y}}$.  In the latter event the two unions
of leaves appearing in (\ref{leaves})
are equal, so in particular, given $\bstar\tilde{x}\in {\sf Lim}_{\hat{x}}$,
there is a unique $\bstar\tilde{y}\in {\sf Lim}_{\hat{y}}$ for which
$L_{\bstar\tilde{x}} =L_{\bstar\tilde{y}}$.  Since $T\cap L_{\hat{x}}$ is
countable, we may thus choose the ordering of each ${\sf Lim}_{\hat{y}}$,
$\hat{y}\in T\cap L_{\hat{x}}$, so that all of the $\alpha$th elements
lie on distinct leaves.  In this way we
may asume that the associated section $T_{\alpha}$ 
intersects any leaf of $\gcom{\mathcal{L}}$ no more than once.
We topologize each section $T_{\alpha}$
through its identification with $T$, and give $\gcom{\mathcal{L}}$
the associated product lamination structure.
By construction of this topology, 
$\bstar p$ becomes an open lamination map.  
\end{proof}

The topology constructed in Theorem~\ref{guctopology} is called
a {\em germ universal cover topology}: it is not unique
and depends on the choice of decomposition (\ref{decomposition}).  
From now on, we assume that
$\gcom{\mathcal{L}}$ has been equipped with such a topology.

There is a canonical simply connected leaf corresponding
to the inclusion $\widetilde{L}\hookrightarrow\gcom{\mathcal{L}}$,
however the other leaves need not be simply connected.  For example, if
$L$ is one of the simply connected leaves of the Sullivan solenoid $\widehat{\D}_{f}$,
then leaves of the associated germ universal cover that correspond to
accumulations of $L$ on an annular leaf will not be simply connected.
Thus $\gcom{\mathcal{L}}$ can be thought of as 
the ordinary universal cover $\widetilde{L}$ surrounded by a nonstandard cloud
of leaves corresponding to the laminar accumulations of $L$; since
these leaves are nowhere dense, one might say that
on passing to $\gcom{\mathcal{L}}$ all of the diophantine 
approximations within $\mathcal{L}$ have been ``unwrapped''.

We now posit $\gcom{\mathcal{L}}$ as the unit space of an enhanced
groupoid structure for $\fgrm (\mathcal{L})$.  
Let $\bast u\in \fgrm (\mathcal{L})$ and $\bstar\tilde{x}\in
\gcom{\mathcal{L}}$.  We say that $\bast u$ {\em acts} on $\bstar\tilde{x}$
if 
there exist representative sequences such that
$\{ u_{\alpha}\cdot \tilde{x}_{\alpha}\}$ defines an $\mathcal{L}$-convergent
sequence $\bast u\cdot\bstar \tilde{x}$ with
\[ \lim (\bast u\cdot\bstar \tilde{x})\;\; =\;\; \lim \bstar \tilde{x}.\]
Defining the domain ${\sf Dom}(\bast u )$ and range 
${\sf Ran}(\bast u )$ of $\bast u$ through this notion
of action, we see that $\gcom{\mathcal{L}}$
yields a new groupoid structure on $\fgrm (\mathcal{L})$, called the
{\em geometric groupoid structure}.  
It is clear that both ${\sf Dom}(\bast u )$ and 
${\sf Ran}(\bast u )$ are sublaminations
of $\gcom{\mathcal{L}}$, since 
$\bstar \tilde{x}\in {\sf Dom}(\bast u )$ implies that $L_{\bstar \tilde{x}}\subset
{\sf Dom}(\bast u )$.  Thus we may view $\fgrm (\mathcal{L})$
as a groupoid of partially defined bijections of $\gcom{\mathcal{L}}$.
Note that the unit space for the old groupoid structure,
$\bstar {\sf D}(\tilde{x},T)$,
maps into the new unit space $\gcom{\mathcal{L}}$ via its
bijection with ${\sf Lim}_{x}$.  There is a canonical
inclusion of the old groupoid structure into the geometric groupoid
structure, given by extension
of domain and range, however in general this map need not be a groupoid
homomorphism.

\begin{assu}  For the remainder of the paper, we will assume
that $\fgrm (\mathcal{L})$ is endowed with the geometric groupoid structure.
\end{assu}

\begin{defi}  We say that $\fgrm (\mathcal{L})$ is {\bf tame}
if whenever $\lim \bstar\tilde{x}=\lim\bstar\tilde{y}$,
there exists $\bast u\in \fgrm (\mathcal{L})$ such that
$\bast u\cdot \bstar\tilde{x}=\bstar\tilde{y}$.
\end{defi}

\begin{prop}  If $\fgrm (\mathcal{L})$ is tame, then the quotient
\[\fgrm ({\cal L})\big\backslash\gcom{\mathcal{L}} \] is
homeomorphic to $\mathcal{L}$.
\end{prop}

\begin{proof}  The equivalence relation enacted by the
action of $\fgrm ({\cal L})$ identifies precisely those points of
$\gcom{\cal L}$ which map to the same point $\hat{x}\in\mathcal{L}$ by
$\bstar p$. Since $\bstar p$ is open, it follows that quotient topology is that the
of $\mathcal{L}$.  
\end{proof}

\begin{theo}\label{actbyhomeo}  If $\fgrm (\mathcal{L})$ is tame and a group, 
then there is a germ universal cover topology on $\gcom{\mathcal{L}}$
for which
$\fgrm (\mathcal{L})$ acts as a group of homeomorphisms.
\end{theo}

\begin{proof}  Let $\gzg{T}$ be the preimage of a transversal $T$
of $\mathcal{L}$.  As $\fgrm (\mathcal{L})$ is a group, 
${\sf Dom}(\bast u )=\gcom{\mathcal{L}}$
for every element $\bast u\in \fgrm (\mathcal{L})$, and moreover
$\bast u (\gzg{T})=\gzg{T}$.  Let $i:T\rightarrow \gzg{T}$ be a section
so that for all $\bstar\tilde{x}\in\gcom{\mathcal{L}}$, 
$i (T)\cap L_{\bstar\tilde{x}}$ contains at most one point.
Since $\fgrm (\mathcal{L})$ acts without fixed points and is tame, 
we have a decomposition as disjoint union
\[   \gzg{T}\;\; =\;\; \bigsqcup_{\bast u\in\fgrm (\mathcal{L})} \bast u (i(T)) .\]
Now construct as in Theorem~\ref{guctopology} a lamination structure
on $\gcom{\mathcal{L}}$ based on this decomposition.  It follows then
that each $\bast u\in \fgrm (\mathcal{L})$ acts homeomorphically
on $\gcom{\mathcal{L}}$.
\end{proof}

\begin{prop}\label{continuity}  Let 
$F:(\mathcal{L}, L)\rightarrow (\mathcal{L}',L')$ be a lamination map, where $L$
and $L'$ are dense leaves.
Then $F$ induces a map
\[ \gcom{F}
:\gcom{\mathcal{L}}\longrightarrow\gcom{\mathcal{L}'},\]
continuous with respect to appropriate choices of germ universal cover topologies.
\end{prop}

\begin{proof}  Denote by $p':\widetilde{L}'\rightarrow L'$ the universal cover.
The map $\gcom{F}$ is defined by representing 
$\bstar\tilde{x}$ by a sequence $\{
\tilde{x}_{\alpha}\}$ and taking $\gcom{F}(\bstar\tilde{x})$ to 
be the asymptotic class of $\{
\widetilde{F}(\tilde{x}_{\alpha})\}$.  Now let $\gzg{\tau'}$ be any germ universal
cover topology on $\gcom{\mathcal{L}'}$, say constructed
from a transversal $T'$.  Since $F$ is a lamination
map, there exists a transversal $T$ with $F(T)\subset T'$.  We
may thus find a decomposition $\gzg{T}=\sqcup T_{\alpha}$
compatible with that of $\gzg{T'}$ {\em i.e.}\  so that
$\gcom{F}(T_{\alpha})\subset T_{\alpha}'$ for all $\alpha$.  Let
$\gzg{\tau}$ to be the associated germ universal cover topology.  Then
$\gcom{F}$ is continuous
with respect to $\gzg{\tau}$ and $\gzg{\tau'}$.
\end{proof}

We now return to the question of functoriality, which we must address
in view of our adoption of a new groupoid structure.  
If we reconsider the notions of
fidelities and trainings with regard to the
geometric groupoid structure, then the analogue of Theorem~\ref{funcstab} --
as well as its corollaries -- remain true with identical proofs.  For the
remainder of the paper, the concepts of fidelity and
training will be understood in the context of the geometric groupoid structure.

The classical universal cover enjoys the property that the lift
$\tilde{f}:\widetilde{X}\rightarrow\widetilde{Y}$ of a map 
$f:X\rightarrow Y$ is $\pi_{1}X$-equivariant.  We now describe conditions
under which the same can be said for a lamination map.  A germ lamination map
$F:\mathcal{L}\rightarrow \mathcal{L}'$ is said to be {\em geometric} if
for all $\bast u\in\fgrm (\mathcal{L})$,
$ \gcom{F}({\sf Dom}(\bast u ))\subset {\sf Dom}\big(\gzg{F}(\bast u )\big)$,
$\gzg{F}:\fgrm (\mathcal{L})\rightarrow\fgrm (\mathcal{L}')$ is a homomorphism
and 
\[ \gcom{F}(\bast u\cdot\bstar\tilde{x})\;\; =\;\;
\gzg{F}(\bast u )\cdot \gcom{F}(\bstar\tilde{x}) .\]
Examples of geometric maps are the projection $\mathcal{L}_{\rho}\rightarrow B$ 
of a suspension onto its base 
and any map of manifolds $f:M\rightarrow M'$.

We say that a lamination $\mathcal{L}$ is {\em geometrically faithful} if
it has a {\em geometric fidelity}: 
a fidelity $\imath :\mathcal{L}\rightarrow X$ which is geometric
and for which $\gcom{\imath}:\gcom{\mathcal{L}}\rightarrow \gcom{X}$
is injective.  In addition $F:\mathcal{L}\rightarrow \mathcal{L}'$ is 
said to be {\em geometrically trained}
if it possesses a training $(\imath ,\imath' ,f)$ where $\imath ,\imath'$ are
geometric fidelities.
For example, the fidelity $\imath: \mathcal{F}_{V}\rightarrow\T^{p+q}$ 
of a linear foliation of a torus is geometric, as well as the projection
of a suspension onto a compact base.

\begin{theo}  Let $F:\mathcal{L}\rightarrow \mathcal{L}'$ be geometrically trained.  
Then $F$ is geometric.
\end{theo}

\begin{proof}  Let $(\imath ,\imath' ,f)$ be a geometric training.  Then
we have
\[ \gcom{\imath}\circ \gcom{F}(\bast u\cdot\bstar \tilde{x})\;\;=\;\;
\gcom{\imath}\big(\gzg{F}(\bast u )\cdot \gcom{F}(\bstar \tilde{x})\big)\]
which implies the result as $\gcom{\imath}$ is injective.
\end{proof}

\begin{coro}  Suppose $F:\mathcal{F}\rightarrow \mathcal{F}'$ is a lamination
map of foliations such that the inclusions into the underlying manifolds
are geometric fidelities.  Then $F$ is geometric.  In particular
any lamination map of linear foliations of torii is geometric.
\end{coro}

\section{Covering Space Theory}

A surjective lamination map
$P:\mathcal{L}\rightarrow\mathcal{L}'$ is called a {\em
lamination covering} if $P|_{L}$ is a covering map for every leaf
$L\subset\mathcal{L}$. A lamination map which is a covering
map in the classical sense is a lamination covering but not all
lamination coverings occur this way {\em e.g.}\
the
projection $\xi :\mathcal{L}\rightarrow B$ of a suspension
onto its base.  We say that $P$ is {\em cover trained} if it has a
training $(\iota ,\iota', p)$ in which $p:X\rightarrow X'$
is a covering map.

\begin{theo}  Let $P:\mathcal{L}\rightarrow\mathcal{L}'$ be a germ lamination
covering that is cover trained.  Then 
\begin{enumerate}
\item The induced map of fundamental germs
\[\gzg{P}:\fgrm (\mathcal{L}) \;\longrightarrow\;\fgrm (\mathcal{L}') \]
is a groupoid monomorphism.
\item  The induced map of germ universal covers
\[ \gcom{P}:\gcom{\mathcal{L}}\;\longrightarrow\;\gcom{\mathcal{L}'}
\]
is an open, injective map with respect to appropriate choices of germ universal
cover topologies. 
\end{enumerate}
\end{theo}

\begin{proof}  The first statement follows
from the definition of training and the fact that $\bast p$ is injective
on $\bast\pi_{1}$.  Let $L$, $L'$ be dense leaves in $\mathcal{L}$,
$\mathcal{L}'$ containing $x$, $x'$.  Then the lift of the
restriction $P|_{L}$, $\widetilde{P}|_{L}
:\widetilde{L}\rightarrow\widetilde{L}'$, is a homeomorphism.
It follows that the induced map $\gcom{P}$ is injective.
$\gcom{P}$ is automatically open with respect to 
the germ universal cover topologies constructed as in 
Proposition~\ref{continuity}.
\end{proof}

\begin{note}  Here is an example when the map $\gcom{P}$ is not surjective.
Take $\mathcal{L}=\R$, $\mathcal{L}'=\mathbb{S}^{1}$ and $P:\R\rightarrow \mathbb{S}^{1}$
the universal cover.  Then $\gcom{\mathcal{L}}\approx \R$ but
$\gcom{\mathcal{L}'}\approx\bstar\R$.
\end{note}

Thus when $P$ is a cover trained, the image
\[\textsf{C}\;\;=\;\;\gzg{P}\Big(\fgrm (\mathcal{L})\Big)\]
is a subgroupoid of $\fgrm (\mathcal{L}')$.  We shall now construct
lamination coverings from subgroups,   
restricting attention to the case where $\fgrm
(\mathcal{L})$ is tame and a group.  Assume that
$\gcom{\mathcal{L}}$ has been given a germ universal cover
topology $\gzg{\tau}$ of the type guaranteed by Theorem~\ref{actbyhomeo}.
Consider a subgroup
$\textsf{C}<\fgrm (\mathcal{L})$ and denote by $\mathcal{L}_{\textsf{C}}$
the quotient $\textsf{C}\backslash\gcom{\mathcal{L}}$.  Note
that $\mathcal{L}_{\textsf{C}}$ decomposes into a disjoint union of leaves.
Let $L_{\textsf{C}}$ be any leaf over the dense leaf $L\subset\mathcal{L}$.
Consider the set $\textsf{X}$ of topologies
$\gzg{\tau_{\mathsf{C}}}$ on $\gcom{\mathcal{L}}$ that satisfy the following
conditions. 
\begin{enumerate}
\item The induced topology $\tau_{\mathsf{C}}$ on $\mathcal{L}_{\textsf{C}}$
defines a (possibly non Hausdorff) lamination structure for which
$\mathcal{L}_{\textsf{C}}$ is dense and $\mathcal{L}_{\textsf{C}}\rightarrow
\mathcal{L}$ is a lamination map.
\item Let $T_{\textsf{C}}$ be any transversal of $\mathcal{L}_{\textsf{C}}$,
and denote by $\gzg{T_{\textsf{C}}}$ its preimage with the induced topology.
If $\textsf{C}\cdot\bstar\tilde{x}$ is contained in $\gzg{T_{\textsf{C}}}$,
then $\textsf{C}\cdot\bstar\tilde{x}$ is not open in the topology of
$\gzg{T_{\textsf{C}}}$.
\item  The identity map $(\gcom{\mathcal{L}},\gzg{\tau_{\mathsf{C}}})\rightarrow
(\gcom{\mathcal{L}},\gzg{\tau})$ is open and $\fgrm (\mathcal{L})$ acts
by homeomorphisms on $(\gcom{\mathcal{L}},\gzg{\tau_{\mathsf{C}}})$.
\end{enumerate}
$\textsf{X}$ is not empty, as it contains
$\gzg{\tau}$.  If we order the elements
of $\textsf{X}$ with respect to inclusion, then $\textsf{X}$ is closed
under chains and so contains a maximal element which we also denote
$\gzg{\tau_{\mathsf{C}}}$,
called a {\em
covering topology}.  Denote by $\tau_{\mathsf{C}}$ the quotient topology
induced by $\gzg{\tau_{\mathsf{C}}}$ on $\mathcal{L}_{\textsf{C}}$.

\begin{theo}\label{intersec}  {\rm $\mathcal{L}_{\textsf{C}}$} 
is Hausdorff with respect to $\tau_{\mathsf{C}}$ and the map 
{\rm $\mathcal{L}_{\textsf{C}}\rightarrow\mathcal{L}$} is a lamination
covering.
\end{theo}
\begin{proof}  Since $\gzg{\tau_{\mathsf{C}}}$ is maximal, for any $\bstar\tilde{x}$,
\[  \textsf{C}\cdot\bstar\tilde{x}\;\; =\;\; 
\bigcap_{\textsf{C}\cdot\bstar\tilde{x}\subset\mathcal{U}\in\gzg{\tau_{\mathsf{C}}}} 
\mathcal{U}.\]
It follows that $\mathcal{L}_{\textsf{C}}$ is Hausdorff.  By construction,
$\mathcal{L}_{\textsf{C}}\rightarrow\mathcal{L}$ is surjective and
a covering when restricted to any leaf.
\end{proof}

Two lamination coverings
$P_{i}:\mathcal{L}_{i}\rightarrow\mathcal{L}$, $i=1,2$, are
{\em isomorphic} if there exists a geometric homeomorphism
$F:\mathcal{L}_{1}\rightarrow \mathcal{L}_{2}$ such that
$P_{1}=P_{2}\circ F$.  The group of automorphisms of a lamination
cover $P$ is denoted ${\sf Aut}(P)$.

\begin{prop}  Let $\bast u\in \fgrm (\mathcal{L},x)$ and {\rm $\textsf{C}'=
\bast u\cdot  \textsf{C}\cdot\bast u^{-1}$}.  Then there exist covering
topologies {\rm $\gzg{\tau_{\mathsf{C}}}$} and {\rm $\gzg{\tau_{\mathsf{C}'}}$} so that
{\rm $\mathcal{L}_{\textsf{C}}\rightarrow\mathcal{L}$} and 
{\rm $\mathcal{L}_{\textsf{C}'}\rightarrow\mathcal{L}$}
are isomorphic.
\end{prop}

\begin{proof}  Choose $\gzg{\tau_{\mathsf{C}}}$ a covering topology for
$\mathcal{L}_{\textsf{C}}$ and let $\gzg{\tau_{\mathsf{C}'}}$ be the image
of $\gzg{\tau_{\mathsf{C}}}$ by $\bast u$.  Then $\gzg{\tau_{\mathsf{C}'}}$
is a covering topology for $\mathcal{L}_{\textsf{C}'}$.  With respect
to these choices, the
bijection 
$\bstar\tilde{x}\;\longmapsto\; \bast u\cdot\bstar\tilde{x}$ defines
a homeomorphism
\[ (\gcom{\mathcal{L}},\gzg{\tau_{\mathsf{C}}})\; \longrightarrow\;
(\gcom{\mathcal{L}},\gzg{\tau_{\mathsf{C}'}})\]
which
descends to a geometric homeomorphism of covers.  
\end{proof}

Now suppose $\textsf{C}\vartriangleleft
\fgrm (\mathcal{L},x)$ is a normal subgroup,
$\gzg{\tau_{\mathsf{C}}}$ the covering topology and
$P_{\textsf{C}}:\mathcal{L}_{\textsf{C}}\rightarrow\mathcal{L}$
the associated covering.

\begin{theo}  {\rm ${\sf Aut}(P_{\textsf{C}})$} is isomorphic to the quotient 
{\rm $\fgrm (\mathcal{L},x)/\textsf{C}$}.  The quotient of 
{\rm $\mathcal{L}_{\textsf{C}}$} by {\rm $\fgrm (\mathcal{L},x)/\textsf{C}$ }
is $\mathcal{L}$.
\end{theo}

\begin{proof}  Every element of $\mathcal{L}_{\textsf{C}}$ is a class 
$\textsf{C}\cdot \bstar\tilde{x}$,
for $\bstar\tilde{x}\in\gcom{\mathcal{L}}$.  The action of 
$\fgrm (\mathcal{L},x)/\textsf{C}$ on such classes is well-defined and yields
a subgroup of ${\sf Aut}(P_{\textsf{C}})$.  On the other hand,
the set $\mathsf{Lim}_{\hat{x}}$ is a $\fgrm(\mathcal{L})$-set
on which any geometric
automorphism acts automorphically.  However the automorphism
group of $\mathsf{Lim}_{\hat{x}}$ is $\fgrm (\mathcal{L},x)/\textsf{C}$, so it
follows that ${\sf Aut}(P_{\textsf{C}})\subset \fgrm (\mathcal{L},x)/\textsf{C}$.
It is clear that the quotient of 
{\rm $\mathcal{L}_{\textsf{C}}$} by {\rm $\fgrm (\mathcal{L},x)/\textsf{C}$ }
is $\mathcal{L}$.
\end{proof}

\bibliographystyle{amsalpha}

\end{document}